\documentclass{article}
\usepackage[utf8]{inputenc}

\usepackage[T1]{fontenc}

\usepackage{mathrsfs}
\usepackage{amsthm,amsmath,amsfonts,amssymb}
\usepackage{graphicx}
\usepackage{enumerate,enumitem}
\usepackage{authblk}
\usepackage{cite}
\usepackage{tikz}
\usepackage{subcaption}

\let\OLDthebibliography\thebibliography
\renewcommand\thebibliography[1]{
  \OLDthebibliography{#1}
  \setlength{\parskip}{0pt}
  \setlength{\itemsep}{0pt plus 0.3ex}
}

\usepackage[activate={true,nocompatibility},final,tracking=true,kerning=true,spacing=true,factor=1100,stretch=15,shrink=15]{microtype}
\microtypecontext{spacing=nonfrench}

\allowdisplaybreaks

\makeatletter
\numberwithin{equation}{section}
\numberwithin{figure}{section}
\theoremstyle{plain}
\newtheorem{thm}{\protect\theoremname}
  \theoremstyle{plain}
  \newtheorem{lemma}[thm]{\protect\lemmaname}
    \newtheorem{prop}[thm]{\protect\propname}

\theoremstyle{definition}
\newtheorem{df}[thm]{Definition}

\theoremstyle{plain}

\numberwithin{thm}{section}

\newtheorem{conj}{Conjecture}

\theoremstyle{remark}
\newtheorem*{rem}{Remark}

\makeatother

\usepackage{babel}
\providecommand{\propname}{Proposition}

\providecommand{\lemmaname}{Lemma}
\providecommand{\theoremname}{Theorem}

\renewcommand{\Im}{\imag}
\renewcommand{\Re}{\real}

\newcommand{\CC}{\mathbb{C}}

\newcommand{\NN}{\mathbb{N}}

\newcommand{\vp}{\varphi}

\let \le \leqslant
\let \leq \leqslant
\let \ge \geqslant
\let \geq \geqslant
\let \epsilon \varepsilon
\let \vp \varphi

\let \setminus \smallsetminus

\usepackage{hyperref}
\hypersetup{
    colorlinks=true, 
    linkcolor=black, 
    urlcolor=black, 
    linktoc=all 
}

\newcommand{\m}[1]{\mathbb{#1}}

\renewcommand\Re{\operatorname{Re}}
\renewcommand\Im{\operatorname{Im}}

 \def \1{\mathbf{1}}

\title{Competitive Hele-Shaw flow and quadratic differentials}
\author{Fredrik Viklund\thanks{KTH Royal Institute of Technology, email: frejo@kth.se } \, and David Witt Nystr\"om\thanks{University of Gothenburg, email: wittnyst@chalmers.se}}
\date{}

\begin{document}

\maketitle
\begin{abstract}
We introduce and investigate a generalization of the Hele-Shaw flow with injection where several droplets compete for space as they try to expand due to internal pressure while still preserving their topology. Droplets are described by their closed non-crossing interface curves in $\mathbb{C}$ or more generally in a Riemann surface of finite type. Our main focus is on stationary solutions which we show correspond to the critical vertical trajectories of a particular quadratic differential with second order poles at the source points. The quadratic differentials that arise in this way have a simple description in terms of their associated half-translation surfaces. Existence of stationary solutions is proved in some generality by solving an extremal problem involving an electrostatic energy functional, generalizing a classic problem studied by Teichm\"uller, Jenkins, Strebel and others. We study several special cases, including stationary Jordan curves on the Riemann sphere. We also introduce a discrete random version of the dynamics closely related to Propp's competitive erosion model, and conjecture that realizations of the lattice model will converge towards a corresponding solution to the competitive Hele-Shaw problem as the mesh size tends to zero.

\end{abstract}
\tableofcontents

\section{Introduction}
\subsection{Background and definition of the model}

\subsubsection*{Classical Hele-Shaw flow}
\begin{figure}

\begin{center}
    
\begin{tikzpicture}

\filldraw[fill=blue!20!white, draw=black] (0,-2) to[out=0,in=-120] (1,-1.5) to[out=60,in=-90] (3,0) to[out=90,in=0] (1,1.5) to[out=180,in=45] (-1,-0.5) to[out=-135,in=90] (-1.25,-1) to[out=-90,in=180] (0,-2) node at (1,0) {$D_t$};

\fill[black] (0,-1) circle (0.05cm) node [anchor=west] { $0$};

\draw[thick,->,blue] (3,0) -- (3.5,0);

\draw[thick,->,blue] (-1,-0.5) -- (-1.6,0.1);

\draw[thick,->,blue] (0,-2) -- (0,-2.8);

\draw[thick,->,blue] (1,1.5) -- (1,2);

\draw[thick,->,blue] (1.6,-2) -- (2.1,-2) node [anchor=west] {$=-\nabla G_{D_t,0}$};

\end{tikzpicture}

\end{center}
\caption{Classical Hele-Shaw dynamics.} \label{fig:1}
\end{figure}
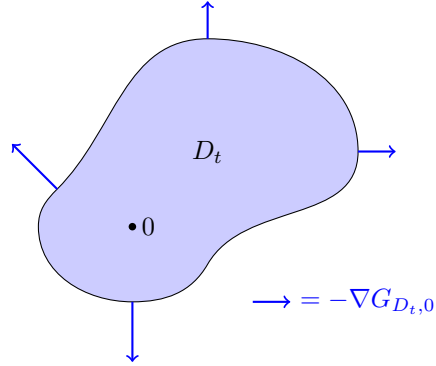

A Hele-Shaw cell consists of two parallel plates separated by a small gap, creating a narrow space into which fluid can be injected or otherwise manipulated. The Hele-Shaw flow is a mathematical model of the propagation of fluid in such a cell, describing the evolution of the fluid interface. The classical setting models a viscous incompressible fluid injected at a constant rate into an otherwise air-filled space through a small hole at the center of the plate. (We do not consider Hele-Shaw flow with suction in this paper.) The region occupied by fluid at time $t$ is identified with a subset $D_t$ of the complex plane $ \mathbb{C} \cong \mathbb{R}^2$, with the point of injection located at the origin and we call $D_t$ a \emph{droplet}. We describe the dynamics of the droplet interface $\partial D_t$. The pressure $p_t$ is assumed to be harmonic in $D_t$ except at the origin, where it has a logarithmic singularity. Hence after normalization we may assume that $-\Delta p_t=\delta_0$, where $\delta_0$ is the Dirac measure at the origin. We furthermore make the simplifying assumption that both the surface tension and the surrounding air pressure is zero, which means that $p_t$ is identically zero on $\partial D_t$. It follows that $p_t=G_{D_t,0}$, where $G_{D_t,0}(z):=G_{D_t}(z,0)$ is the Green's function of $D_t$ with Dirichlet boundary condition and singularity at $0$. 
We say that a family of domains $D_t$ in $\mathbb{C}$ containing the origin is a (classical) solution to the Hele-Shaw problem on the interval $I$ if for all $t \in I$, the normal velocity of $\partial D_t$ at any regular point $z$ is equal to $-\nabla G_{D_t,0}(z)$. See Figure~\ref{fig:1}. Note that $|\nabla G_{D_t,0}(z)|$ is the value of the Poisson kernel for $D_t$ with interior point $0$ and boundary point $z$. If $\partial D_0$ is an analytic Jordan curve, a solution $D_t$ to the Hele-Shaw problem starting from $D_0$ exists at least for small time (see below). For a thorough treatment of the classical Hele-Shaw flow with many references we refer the reader to the book \cite{GV}.

\subsubsection*{Competitive Hele-Shaw flow}

We now introduce a variant of the Hele-Shaw flow where several droplets compete for space  while preserving their topology: the \emph{competitive Hele-Shaw flow}. We first consider the model on the complex plane $\CC$. Let $\gamma$ be a piecewise regular\footnote{A (smooth) regular curve is a $C^\infty$-smooth immersion of an interval viewed as a manifold with boundary, and a piecewise regular curve is a concatenation of finitely many regular curves.} simple closed curve in $\CC$, or more generally a piecewise regular non self-crossing\footnote{A non self-crossing curve is a curve of which there exist arbitrarily small deformations that make it simple.} closed curve in $\mathring{\CC}$, where $\mathring{\CC}$ denotes the differentiable surface with boundary one gets by adding a circle at infinity to $\CC$. (See Section \ref{subsec:prel} for the precise definition.) The reason for letting $\gamma$ be a curve in $\mathring{\CC}$ rather than in $\CC$ or $\hat{\CC}:=\CC\cup\{\infty\}$ is that we want to allow $\gamma$ to pass through infinity but, for reasons that will become apparent later, it should not be possible to continuously deform $\gamma$ across infinity. The interior of $\gamma$ is then a disjoint union $D=\cup E_j$ of simply connected domains in $\CC$. We will call $D$ a \emph{droplet} and the curve $\gamma$ a \emph{droplet interface}.\footnote{While it would be interesting to consider the model with other droplet topologies allowed, here we will restrict ourselves to the simplest setting.} The Green's function $G_{D,w}(z)$ is naturally defined so that if $z$ and $w$ lie in the same simply connected component $E$ of $D$ then $G_{D,w}(z):=G_{E,w}(z)$, while if they lie in different components, $G_{D,w}(z):=0$ by definition.

Let $D$ be a droplet and let a finite number of source points $z_j$, $j=1,\ldots,m$ in $D$, with associated source strengths (weights) $a_j>0$ be given. It is convenient to let $d$ denote the associated source divisor $d:=\sum_{j=1}^m a_j z_j$, and to write $$G_{D,d}(z):=\sum_{j=1}^m a_j G_{D,z_j}(z).$$ By the support $\textrm{supp}(d)$ of $d$ we mean the set of points $\{z_j: j=1,\ldots, m\}$. Note that if $D$ has several components $E_j$, then if $z \in E_j$, $G_{D,d}(z) = G_{E_j,d_j}(z)$ where $d_j$ is the restriction of $d$ to $E_j$.

\begin{df}[Admissible droplet configuration]
A collection $(D^i,\gamma^i,d^i)_{i=1}^n$ of droplets $D^i$ with droplet interfaces $\gamma^i$ and source divisors $d^i$ is said to be \emph{admissible} if the droplets are all disjoint, the droplet interfaces do not cross each other, and for each $i$ we have that $\textrm{supp}(d^i)\subset D^i$. 
\end{df} 

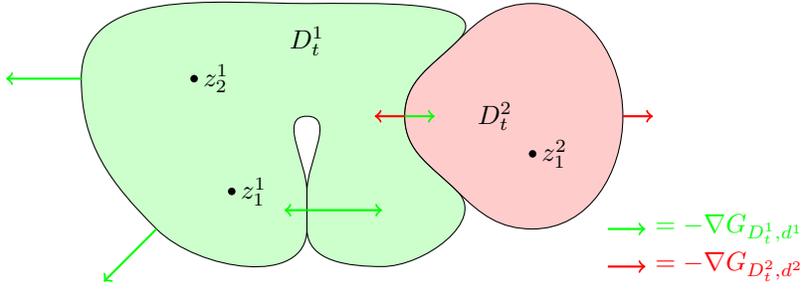
\begin{figure}

\begin{center}
    
\begin{tikzpicture}

\filldraw[fill=green!20!white, draw=black] (0,-0.5) to[out=180,in=90] (0,-1.5) to[out=-90, in=90] (0,-2)[out=-90, in=90] (0,-2) to[out=-90,in=-45] (-2,-2) to[out=135,in=-90] (-3,0) to[out=90,in=180] (0,1)  to[out=0,in=45] (2,0.5) to[out=-135,in=90] (1.3,-0.5) to[out=-90,in=135] (2,-1.5) to[out=-45,in=0] (1,-2.5) to[out=180,in=-90] (0,-2) to[out=90, in=-90] (0,-1.5) to[out=90,in=0] (0,-0.5);

\filldraw[fill=red!20!white, draw=black] (2,0.5) to[out=-135,in=90] (1.3,-0.5) to[out=-90,in=135] (2,-1.5) to [out=-45,in=180] (3,-2) to[out=0,in=-90] (4.2,-0.5) to[out=90,in=0] (3,1) to[out=180,in=45] (2,0.5);

\fill[black] (-1,-1.5) circle (0.05cm) node [anchor=west] { $z^1_1$};

\fill[black] (-1.5,0) circle (0.05cm) node [anchor=west] { $z^1_2$};

\fill[black] (3,-1) circle (0.05cm) node [anchor=west] { $z^2_1$};

\path (0,0.5) node {$D^1_t$} (2.5,-0.5) node {$D^2_t$};

\draw[thick,->,green] (0,-1.75) -- (1,-1.75);
\draw[thick,->,green] (0,-1.75) -- (-0.3,-1.75);
\draw[thick,->,green] (-2,-2) -- (-2.7,-2.7);
\draw[thick,->,green] (-3,0) -- (-4,0);
\draw[thick,->,red] (1.3,-0.5) -- (0.9,-0.5);
\draw[thick,->,green] (1.3,-0.5) -- (1.7,-0.5);
\draw[thick,->,red] (4.2,-0.5) -- (4.6,-0.5);

\draw[thick,->,red] (4,-2.5) -- (4.5,-2.5) node [anchor=west] {$=-\nabla G_{D^2_t,d^2}$};
\draw[thick,->,green] (4,-2) -- (4.5,-2) node [anchor=west] {$=-\nabla G_{D^1_t,d^1}$};

\end{tikzpicture}

\end{center}
\caption{Competitive Hele-Shaw dynamics.} \label{fig:2}
\end{figure}

We will now describe the competitive Hele-Shaw flow on droplet configurations. Roughly speaking, each droplet tries to expand according to the Hele-Shaw dynamics with injection at the source points (with weights) but if two droplets meet, they push against each other and so ``compete'' to expand. If one droplet lies on both sides of an interface, the droplet pushes against itself. See Figure~\ref{fig:2}. If two interfaces do intersect they will continue to move together, ensuring that the droplets continue to be disjoint and the interfaces continue to not cross. In this sense, the dynamics preserve the topology of the system of droplets. 

To give the precise definition of the competitive Hele-Shaw flow we need some auxiliary notions.
\begin{itemize}
\item{
We say that a point $z\in \CC$ is a \emph{regular interface point} of an admissible droplet configuration $(D^i,\gamma^i,d^i)_{i=1}^n$ if it belongs to at least one droplet interface and it is a regular point of each droplet interface it belongs to. If $z$ is a regular interface point and $z$ lies in the closure of a droplet $D^i$, then we say that $D^i$ is a \emph{neighboring droplet} of $z$. If $D^i$ is a neighboring droplet of $z$, expressions such as $\nabla G_{D^i,d^i}(z)$ are understood in the sense of non-tangential limits at (the prime end) $z$ taken from $D^i$. Note that a regular interface point can belong to many droplet interfaces, but it can have at most two neighboring droplets. }

\item{We say that a family of closed piecewise regular curves $\gamma_t$, $t\in I$, is \emph{differentiable} if there is a continuous function $(I \times S^1) \ni (t,\theta) \mapsto \gamma_t(\theta)$ such that for each $t$, $\theta \mapsto \gamma_t(\theta)$ is a piecewise regular parametrization and for each regular point $\theta$, $t \mapsto  \gamma_t(\theta)$ is differentiable. }

\item{Let $\gamma_t$, $t\in I$ be a differentiable family of closed piecewise regular curves and let $z=\gamma_{t_0}(\theta)$ be a regular point of $\gamma_{t_0}$. The \emph{normal velocity} of $\gamma_{t_0}$ at $z$ is then defined as the normal component of $\dot{\gamma}_{t_0}(\theta)$.}
\end{itemize}

\begin{figure}

\begin{center}
    
\begin{tikzpicture}

\filldraw[fill=blue!20!white, draw=black] (0,-0.25) to[out=45,in=180] (1,0) to[out=0,in=0] (1,2) to[out=180,in=90] (-2,0) to[out=-90,in=180] (1,-2) to[out=0,in=0] (1,0) to[out=180,in=45] (0,-0.25) node at (0,0.5) {$D^i_t$};

\fill[red] (0,-0.25) circle (0.05cm);

\filldraw[fill=green!20!white, draw=black] (5,0) to[out=0,in=0] (5,2) to[out=180,in=90] (2,0) to[out=-90,in=180] (5,-2) to[out=0,in=0] (5,0) node at (4,0.5) {$\tilde{D}^i_t$};

\end{tikzpicture}

\end{center}
\caption{Removal of slit singularity.} \label{fig:3}
\end{figure}
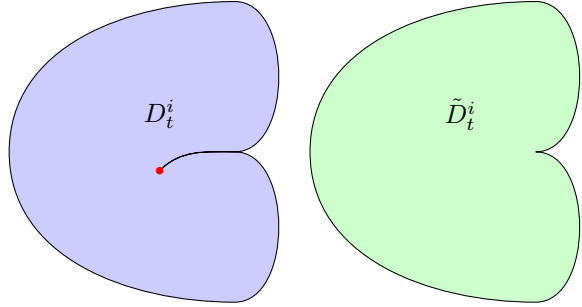

\begin{df}[Competitive Hele-Shaw problem]\label{def:CHS-dynamics} We say that a differentiable family of admissible droplet configurations $(D^i_t,\gamma^i_t, d^i)_{i=1}^n, \, t \in I,$ is a \emph{solution to the competitive Hele-Shaw problem with source divisors $d^i$} if for all $i$ and $t \in I$ the normal velocity of the interface $\gamma^i_t$ at any regular interface point $z \in \gamma^i_t \cap \CC$ is given by \begin{align}\label{CHS1}
V^i_t(z)=-\sum\nabla G_{D^j_t,d^j}(z),\end{align} where the sum is taken over the neighboring droplets $D^j$ of $z$. In particular, if there are no neighboring droplets we have that $V^i_t(z)=0$. If there is one neighboring droplet $D^j$ of $z$ that lies on both sides of the interface (i.e., $z$ corresponds to two distinct prime ends of $D^j$), the sum in \eqref{CHS1} is interpreted as the sum of the gradients taken from the two different sides of the interface, in the sense of non-tangential limits. 
\end{df}
As for the classical Hele-Shaw problem, we have only specified the dynamics at regular points. However, one singular case requires special attention. Consider a droplet interface $\gamma^i_t$ with a singularity as shown in Figure \ref{fig:3}. At the tip of the slit (marked with a red dot in Figure \ref{fig:3}) we have $|\nabla G_{D^i_t,d^i}|=\infty$, which intuitively should mean that the self-intersecting part of the interface retracts with infinite speed. Therefore we replace $D^i_t$ with $\tilde{D}^i_t$. However, if $\infty$ lies on the slit part of $\gamma^i_t$, then, since $\gamma^i_t$ cannot cross $\infty$, we only retract the slit up to $\infty$. See also the second remark after Definition~\ref{def:CHS-surface}. We add the following to Definition~\ref{def:CHS-dynamics}:
\begin{itemize}
\item{If for some $i$ and $t$ the interface $\gamma^i_t$ has an inward-pointing slit singularity as in Figure \ref{fig:3} and $\tilde{D}^i_t$ is the droplet one gets after removing that slit singularity (possibly up to the point at infinity), then $\lim_{s\to t+}D^i_s=\tilde{D}^i_t$.}
\end{itemize}
\begin{rem}
Definition~\ref{def:CHS-dynamics} can be naturally extended to other singular points assuming the gradient limits taken from all adjacent droplet components exist and are finite (i.e., the limits at all prime ends exist and are finite). This happens, e.g., if $k > 2$ interfaces meet at the point and the intersection angles are all acute -- in this case all gradients are $0$. Such configurations, with all angles equal, appear generically in stationary solutions. 
\end{rem}

We shall comment on existence of the competitive Hele-Shaw flow below. In this paper our main focus on is stationary solutions to the competitive Hele-Shaw problem, i.e., droplet configurations invariant with respect to the dynamics in Definition~\ref{def:CHS-dynamics}. 
\begin{df}[Stationary solution in the plane]\label{def:stationary-solution0} We say that $(D^i,\gamma^i)$ is a stationary solution to the competitive Hele-Shaw problem with source divisors $d^i$ if  $(D^i,\gamma^i,d^i)$ is an admissible droplet configuration such that at any regular interface point $z \in \gamma^i\cap \CC$ the vector field $V^i(z)$ as in \eqref{CHS1} equals $0$, and no droplet interface has a slit singularity as in Figure \ref{fig:3} unless the tip is at $\infty$.
\end{df}

\subsubsection*{Competitive Hele-Shaw flow on a Riemann surface}

Let $\Sigma$ be a Riemann surface of finite type, i.e., a closed Riemann surface with no or finitely many punctures. We always equip $\Sigma$ with a complete Riemannian metric $g$ compatible with the complex structure and we write $(\Sigma, g)$ for this Riemannian surface.  

If $\Sigma$ is non-compact, and hence punctured, let $\tilde{\Sigma}$ denote the associated closed unpunctured Riemann surface. We also let $\mathring{\Sigma}$  denote the differentiable surface with boundary obtained by adding a circle to $\Sigma$ at each puncture $p_i$, see Section \ref{subsec:prel}. 

Let $\gamma$ be a closed null-homotopic piecewise regular non self-crossing curve in $\Sigma$ (or in $\mathring{\Sigma}$ if $\Sigma$ is non-compact). In  both cases, the interior $\textrm{int}(\gamma)$ of $\gamma$ is a disjoint union of simply connected domains in $\Sigma$. (If there are two possible choices of interiors for a given $\gamma$ we pick the one which makes the orientation positive.) We will call such a set a \emph{droplet} and the curve $\gamma$ a droplet interface. The Green's function for a droplet (possibly with several components) is defined analogously to the case of the complex plane. The gradient $\nabla_g G_{D,d}$ now depends on the metric but having chosen this we may define a vector field exactly as in \eqref{CHS1}.

As before we say that a droplet configuration $(D^i,\gamma^i,d^i)_{i=1}^n$ is \emph{admissible} if the droplets are all disjoint, the droplet interfaces do not cross each other, and for each $i$ we have that $\textrm{supp}(d^i)\subset D^i$. 

Now we can define the competitive Hele-Shaw problem on $(\Sigma,g)$ just as in $\CC$ with the Euclidean metric.

\begin{df}[Competitive Hele-Shaw flow on a surface]\label{def:CHS-surface}
We say that a differentiable family of admissible droplet configurations $(D^i_t,\gamma^i_t, d^i)_{i=1}^n, \, t \in I,$ in $\Sigma$ is a \emph{solution to the competitive Hele-Shaw problem with source divisors $d^i$} if for all $i$ and $t \in I$ the normal velocity of the interface at any regular interface point $z \in \gamma^i_t\cap \Sigma$ is given by the vector field as in \eqref{CHS1} with $\nabla$ replaced by $\nabla_g$. Furthermore, if for some $i$ and $t$ the interface $\gamma^i_t$ has a singularity as in Figure \ref{fig:3} and $\tilde{D}^i_t$ is the droplet one gets after removing the slit (possibly up to a puncture), we require that $\lim_{s\to t+}D^i_s=\tilde{D}^i_t$.
\end{df}
\begin{rem}
    In general a solution will depend on the choice of metric $g$. The physical meaning of the choice of  metric can be understood locally in terms of a varying permeability, as was first considered in the classical Hele-Shaw setting by Hedenmalm-Shimorin in \cite{HS} (see also \cite{RWN}).
\end{rem}
\begin{rem}
The effect of introducing punctures is that the droplet interfaces $\gamma^i_t$ are prevented from passing through the punctures, so a puncture will thus act as a fixed obstacle. In particular, if $\gamma^i_t$ develops a slit singularity as in Figure \ref{fig:3} where the tip of slit (marked with red) lies at the puncture we do not replace $D^i_t$ by $\tilde{D}^i_t$, as this would require $\gamma^i_t$ to cross the puncture. If $\gamma^i_t$ develops a singularity as in Figure \ref{fig:3} and there is a puncture somewhere in the middle of the self-intersecting part, then we remove the self-intersecting part of the interface up to that puncture. 
\end{rem}
\begin{df}[Stationary solution on a surface]\label{def:stationary-solution} We say that $(D^i,\gamma^i)$ is a stationary solution to the competitive Hele-Shaw problem with source divisors $d^i$ on $\Sigma$  if $(D^i,\gamma^i,d^i)$ is an admissible  droplet configuration such that at any regular interface point on $\gamma^i\cap \Sigma$ the vector field defined as in \eqref{CHS1} with $\nabla$ replaced by $\nabla_g$ equals $0$, and no interface has a slit singularity as in Figure \ref{fig:3} unless the tip is at a puncture. 
\end{df}

Stationary solutions, in contrast to non-stationary solutions, do not depend on the metric $g$. Indeed, if a sum of the form $-\sum \nabla_g G_{D^j_t,d^j}=0$ for some $g$, it equals $0$ for any $g$. Therefore we usually do not specify $g$ when considering a stationary solution.

\subsection{Existence of the competitive Hele-Shaw flow}

Even for the classical Hele-Shaw problem (with injection) with analytic starting data, short-time existence is quite non-trivial. When the initial domain $D_0$ is simply connected with a real-analytic and non self-intersecting boundary, short-time existence was first obtained by Kufarev-Vinogradov \cite{VK}. A more modern and simpler proof based on the abstract Cauchy-Kovalevskaya theorem was later given by Reissig-von Wolfersdorf \cite{RvW}. Short-time existence in the case of a general metric $\textit{g}$ was established in \cite{RWN}. 

For the classical Hele-Shaw problem there exists a notion of weak solution and existence and uniqueness of weak solutions have been established in a very general setting, moreover, the weak solution is a ``strong'' solution as long as it is smooth, see \cite{GV}. We do not know how to formulate a useful notion of weak solution for the competitive Hele-Shaw problem and this remains an interesting open problem. We expect that one can obtain local existence in the case of two droplets in $\hat{\mathbb{C}}$ separated by an analytic Jordan curve along similar lines as in \cite{RvW}.

It seems reasonable to expect the competitive Hele-Shaw flow to locally be at least as smoothing as the classical Hele-Shaw flow. Note that after removing a slit singularity as in Figure \ref{fig:3} a cusp-like singularity may remain. Explicit examples of the classical Hele-Shaw flow such as the Polubarinova-Galin cardioid (see \cite{GV}) suggest that such inward-pointing cusp singularities are instantly resolved. We expect this to be the case here as well and that in fact solutions exist for all time and converge towards stationary solutions.

\begin{conj} \label{conj:cHS}
    Let $(D^i_0,\gamma^i_0,d^i)_{i=1}^n$ be an admissible droplet configuration in $\Sigma$ (if $\Sigma=\hat{\CC}$ we suppose that $n>1$). Then there exists a unique solution $(D_t^i, \gamma^i_t)$, $t\in [0,\infty)$, to the competitive Hele-Shaw problem on $(\Sigma,g)$ with source divisors $d^i$ starting from $(D_0^i,\gamma_0^i)$. Furthermore (if $\Sigma=\CC$ we suppose that $n>1$ while if $\Sigma=\hat{\CC}$ we suppose that $n>2$ ), the solution $(D_t^i, \gamma^i_t)$ will converge as $t \to \infty$ to a stationary solution $(D^i, \gamma^i)$ to the competitive Hele-Shaw problem on $(\Sigma,g)$ with source divisors $d^i$.   
\end{conj}
In Section~\ref{sect:introCIE} we introduce a random lattice version of the competitive Hele-Shaw flow which can be viewed as a probabilistic regularization. It has the advantage of being immediately well-defined and some experimental support for the conjecture will be given when considering simulations of this discrete model, see Sections \ref{sect:introCIE} and \ref{sect:discretemodels}.

\subsection{Stationary solutions and quadratic differentials} \label{subsec:st}
In this paper we focus on stationary solutions to the competitive Hele-Shaw problem. We will establish and study in particular a link to quadratic differentials. This link is most easily seen assuming existence of a stationary solution. 
\subsubsection*{Quadratic differential from stationary solution} Let $\Sigma$  be a Riemann surface of finite type. Suppose $(D^i,\gamma^i)$ is a stationary solution to the competitive Hele-Shaw problem on $\Sigma$ with source divisors $d^i, i=1,\ldots, n$. Then for each $i$, $\partial_z G_{D^i,d^i}(z)dz$ is a well-defined meromorphic abelian differential on each droplet component. It has poles of order one at the source points with residues given by $-1/2$ times the source strengths. For any regular point $z \in \partial D^i \cap \partial D^j$ on an interface, the stationarity condition implies \[\partial_z G_{D^i,d^i}(z)dz=-\partial_z G_{D^j,d^j}(z)dz.\] Given this observation, the following result is not hard to prove. 
\begin{prop}\label{proB}
Define locally in each droplet component of $D^i$ a quadratic differential
 $$\vp(z)dz^2=\left(\partial_z G_{D^i,d^i}(z)\right)^2 dz^2.$$
 Then $\vp$ extends to a meromorphic quadratic differential on the whole surface $\Sigma$ (or on the closed surface $\tilde{\Sigma}$ if $\Sigma$ is non-compact).
\end{prop}
Proposition~\ref{proB} is proved in Section~\ref{sect:stationary-to-QD}

\subsubsection*{Stationary solution from quadratic differential} The result in the previous section assumed the existence of a stationary solution, but it is not a priori clear that any stationary solutions exist.
Our main result is the following.
\begin{thm}[Existence of stationary solutions] \label{ThmA}
   Let $\Sigma$ be a Riemann surface of finite type. Let $(D^i_0,\gamma^i_0,d^i)_{i=1}^n$ be an admissible droplet configuration in $\Sigma$. If $\Sigma=\CC$ suppose that $n>1$ and if $\Sigma=\hat{\CC}$ suppose that $n>2$. Then there exists a stationary solution $(D^i, \gamma^i)_{i=1}^n$ to the competitive Hele-Shaw problem on $(\Sigma,g)$ with source divisors $(d^i)_{i=1}^n$ such that for all $i$, $\gamma^i$ is homotopic to $\gamma^i_0$ in $\Sigma\smallsetminus \cup_i \textrm{supp}(d^i)$ (or in $\mathring{\Sigma}\smallsetminus \cup_i \textrm{supp}(d^i)$ if $\Sigma$ is non-compact).
\end{thm}

The proof of Theorem \ref{ThmA} is given in Section~\ref{sect:existenceproof}. The main step in the proof is to construct a meromorphic quadratic differential related to the competitive Hele-Shaw problem as above. Existence of this quadratic differential is obtained by solving an extremal problem involving an electrostatic energy functional defined as follows.

Let $D$ be a droplet in $\Sigma$ and $d = \sum_{j=1}^n a_j z_j$ a divisor with support in $D$. Choose a local coordinate near each marked point. We define the \emph{reduced Green's energy} of $(D,d)$ by
\begin{align}\label{def:Green-Energy0}
I(D, d) := \sum_{1 \le j \neq k \le n} a_ja_kG_D(z_j,z_k) + \sum_{j=1}^n a_j^2M_{D}(z_j),
\end{align}
where $M_D(z_j)$ is the reduced modulus of $D\smallsetminus\{z_j\}$ which depends on the choice of coordinate. (See Section \ref{subsec:prel} for the definition of reduced modulus.) If $\mathcal{D} = (D^i,\gamma^i,d^i)_{i=1}^n$ is an admissible droplet configuration we define its reduced Green's energy by $$I(\mathcal{D} ):=\sum_{i=1}^n I(D^i,d^i).$$

If $\mathcal{D}_0 = (D_0^i,\gamma_0^i,d^i)_{i=1}^n$ is an admissible droplet configuration we let $\mathcal{H}(\mathcal{D}_0)$ denote the space of all admissible droplet configurations $(D^i,\gamma^i,d^i)$ such that for each $i$, $\gamma^i$ is homotopic to $\gamma^i_0$ in $\Sigma\smallsetminus \cup_i \textrm{supp}(d^i)$ (or in $\mathring{\Sigma}\smallsetminus \cup_i \textrm{supp}(d^i)$ if $\Sigma$ is non-compact). 
\begin{prop}[Solution of extremal problem]\label{thm2}
There exists a $\mathcal{D}_{\infty}\in \mathcal{H}(\mathcal{D}_0)$ which maximizes the reduced Green's energy in $\mathcal{H}(\mathcal{D}_0)$. The quadratic differential which in each droplet component of $\mathcal{D}_{\infty}$ is defined by \[
    \vp(z)dz^2 = (\partial_z G_{D^{i}, d^i}(z))^2 dz^2
    \]
is almost everywhere equal to a meromorphic quadratic differential on the whole surface $\Sigma$. The vertical trajectories connecting the finite critical points of $\varphi$ trace the droplet interfaces of $\mathcal{D}_{\infty}$. 
\end{prop}
The finite critical points are the zeros and the first order poles. Proposition~\ref{thm2} is proved in Section~\ref{sect:existenceproof}. The proof is based on compactness and quasiconformal variation combined with Weyl's lemma but proving that the quadratic differential corresponds to an admissible droplet configuration in $\mathcal{H}(\mathcal{D}_0)$ also requires a topological argument. Given this, the proof of Theorem~\ref{ThmA} is almost immediate: The vertical foliation of $\vp$ is given by the level sets of the Green's functions $G_{D^i,d^i}$,
 \begin{eqnarray*}
 \textrm{Re}\left(\int_{z_0}^z \sqrt{\vp(w)}dw\right)=\int_{z_0}^z \textrm{Re}\left(\partial_w G_{D^i,d^i}dw\right)=\int_{z_0}^z \frac{1}{2}dG_{D^i,d^i}=\frac{1}{2}(G_{D^i,d^i}(z)-G_{D^i,d^i}(z_0)).
 \end{eqnarray*}

\begin{rem}When all source divisors are singletons the reduced Green's energy is simply a weighted sum of reduced moduli. This exact functional was used by Teichm\"uller in the case of two punctured discs in the plane (see the discussion on p99 of \cite{strebel-book}) and later by Jenkins to prove a version of Theorem~\ref{thm2} in this special case, see Remark~1 in \cite{Jenkins57}. In this case solutions are unique. However, the more general setting we consider here and the reduced Green's energy seem to be new and we have not proved uniqueness.
\end{rem}

\begin{figure}

\begin{center}
    
\begin{tikzpicture}

\fill[green!20!white] (-3,2) -- (6,2) -- (6,-2) -- (-3,-2) -- (-3,2);

\draw (4,0) -- (6,0);
\path (-1,0.5) node {$D^0$};

\filldraw[fill=blue!20!white, draw=black] (4,0) -- (3.97 ,0.05196152) -- (3.93871865,0.09679904) -- (3.89199263,0.15657086) -- (3.82298694,0.23139513) --(3.71985828, 0.3200129) -- (3.66249479,  0.359608403) -- (3.63205849, 0.378174918) -- (3.59877223, 0.396674371) -- (3.56227639, 0.414894507) -- (3.56025896,0.41584071) -- (3.52214906,  0.432549765) -- (3.47789333, 0.449251957) -- (3.42892285,  0.464466534) -- (3.37454713, 0.477446835) -- (3.31396092, 0.487131592) -- (3.29687236,0.48895403) --
(3.24624957, 0.491983074) -- (3.17044146, 0.489726886) --
(3.08568697, 0.476933411) -- (2.99176776, 0.448372611) -- (2.89043979, 0.396192296) -- (2.84151581,0.36022014) --(2.78863837, 0.309720878) -- (2.70414728, 0.179431435) -- (2.66673832, 0.114869060) -- (2.655, 0) --(2.66673832, -0.114869060) --
(2.69595722, -0.159306980) -- (2.77615632, -0.295405601) -- (2.84440779,-0.3625962) --
(2.87701595, -0.387158061) -- (2.97896308, -0.443138290) --
(3.07401008, -0.474255723) -- (3.15996475, -0.488752190) -- (3.23689186, -0.492151046) -- (3.29872478,-0.48877607) --(3.30559882, -0.488073707) --
(3.36705578, -0.478916620) -- (3.42218929, -0.466295735) --
(3.47181968, -0.451323364) -- (3.51665193, -0.434780253) --
(3.55728506, -0.417224069)   -- (3.56133943,-0.41533473) -- (3.7205398,-0.31950391) -- (3.82343732,-0.230951) -- (3.89229674,-0.15620666) -- (3.93892377,-0.09651615) -- (4,0) node at (3.4,0.2) {$D^1$};

\fill[black] (3,0) circle (0.05cm) node [anchor= north west] { $1$};
\fill[black] (0,0) circle (0.05cm) node [anchor=north west] { $0$};

\end{tikzpicture}

\end{center}
\caption{An explicit example of a stationary solution with two droplets in $\CC$. The source point at $0$ has weight $2$ while the source point at $1$ has weight $1$. The self-intersecting part of the interface extends to the puncture at $\infty$. This example corresponds to an extremal problem and quadratic differential considered by Teichm\"uller.} \label{fig:4}
\end{figure}
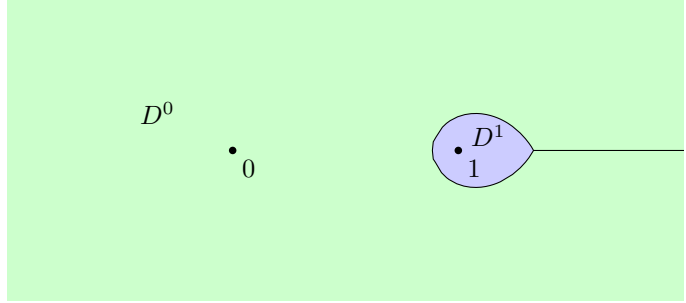

The link to quadratic differentials leads to explicit descriptions of some examples of stationary solutions to the competitive Hele-Shaw problem. For instance the quadratic differential \begin{align}\label{example:teichmuller}\varphi(z)dz^2=\frac{4-3z}{z^2(z-1)^2}dz^2\end{align} on $\CC$ can be shown to correspond to a stationary solution $(D^i,\gamma^i)$ with source divisors $d^1=2 \cdot 0 $ and $d^2=1 \cdot 1$, and from this the droplet interfaces can easily be computed (see Figure \ref{fig:4}). In particular the quadratic differential is seen to have a zero at $z=4/3$, hence that is the point where $D^0$ meets $D^1$ as well as itself. We believe the interface erosion droplet configuration in the simulation in Figure~\ref{fig:sim} approaches the trajectories of this particular quadratic differential. More details of these computations, and more explicit examples, are given in Section \ref{sec:ex}.

\subsubsection*{Stationary solutions and half-translation surfaces}
It is well-known that a quadratic differential on $\Sigma$ gives rise to a representation of $\Sigma$ as a half-translation surface, i.e., a surface given as a collection of polygons with the edges identified by translations or minus-translations (such as $z\mapsto -z+b$). In fact, these two descriptions are equivalent, see, e.g., \cite{Zor}.

\begin{figure}

\begin{center}
    
\begin{tikzpicture}

\fill[blue!20!white] (-5,1) -- (0,1) -- (0,0) -- (-5,0) -- (-5,1);

\path (-2.5,0.5) node {$S^i$};

\fill[green!20!white] (-5,-0.5) -- (0,-0.5) -- (0,-2) -- (-5,-2) -- (-5,-0.5);

\path (-2.5,-1.25) node {$S^j$};

\draw[blue] (-5,1) -- (0,1);

\draw[blue] (-5,0) -- (0,0);

\draw[green] (-5,-0.5) -- (0,-0.5);

\draw[green] (-5,-2) -- (0,-2);

\draw (0,1) -- (0,0.7);

\draw[red, dashed] (0,0.7) -- (0,0.3);

\draw (0,0.3) -- (0,0);

\draw (0,-0.5) -- (0,-1);

\draw[red, dashed] (0,-1) -- (0,-1.4);

\draw (0,-1.4) -- (0,-2);

\fill[blue!20!white] (1.5,0.5) to[out=-120,in=120] (1,-1.5) to[out=0,in=180] (3,-1) to[out=60,in=-60] (3,0.7) to[out=180,in=0] (1.5,0.5);

\fill[black] (2.5,-0.2) circle (0.05cm) node [anchor=west] { $z_i$};

\path (2,-0.5) node {$D^i$};

\fill[green!20!white] (3,0.7) to[out=60,in=150] (5.5,0.5) to[out=-90,in=60] (4.5,-1.5) to[out=180,in=-60] (3,-1) to[out=60,in=-60] (3,0.7);

\fill[black] (4,-0.75) circle (0.05cm) node [anchor=west] { $z_j$};

\path (4.2,0) node {$D^j$};

\draw[red, dashed] (3,-1) to[out=60,in=-60] (3,0.7);

\draw (3,0.7) to[out=180,in=0] (1.5,0.5) to[out=-120,in=120] (1,-1.5) to[out=0,in=180] (3,-1);

\draw (3,0.7) to[out=60,in=150] (5.5,0.5) to[out=-90,in=60] (4.5,-1.5) to[out=180,in=-60] (3,-1);

\draw (1.5,0.5) to[out=120, in=-60] (1.3,0.9);

\draw (1,-1.5) to[out=-120, in=60] (0.8,-1.9);

\draw (1,-1.5) to[out=-120, in=60] (0.8,-1.9);

\draw (4.5,-1.5) to[out=-60,in=120] (4.7,-1.9);

\draw (5.5,0.5) to[out=30,in=-150] (5.9,0.7);

\end{tikzpicture}

\end{center}
\caption{The associated half-translation surface of a stationary solution.} 
\label{fig:5}
\end{figure}
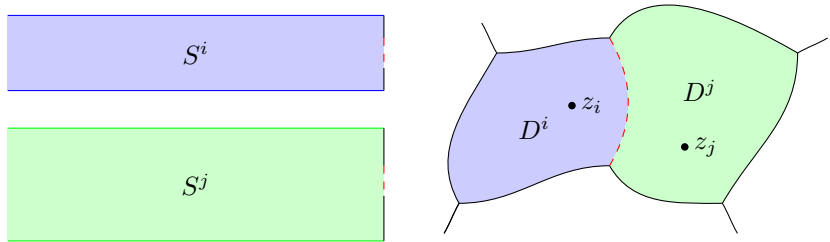

Let $\vp$ be the quadratic differential associated to a stationary solution $(D^i,\gamma^i)$ where each source divisor $d^i$ is a singleton divisor $d^i=a^iz^i$. As is shown in Section \ref{subsec:half}, $\vp$ yields a representation of $\Sigma$ where $D^i$ corresponds to the half-strip $S^i:=\{\Re(z)<0, 0\leq \Im(z)\leq \pi a^i\}$ with the top and bottom of each half strip being identified by the translation $z\mapsto z+\pi ia^i$. Furthermore, if the interface $\gamma^i$ intersects the interface $\gamma^j$ along an arc, then this corresponds to a piece of the vertical boundary of $S^i$ being identified via a minus-translation to a piece of the vertical boundary of $S^j$, see Figure \ref{fig:5}. 

When the divisors are not all singletons a similar description holds, but instead of only half-strips there are now some extra rectangles glued in. We characterize  such surfaces -- half-translation surfaces of \emph{Green's type}. 
\begin{thm}\label{thm:1-1}
Stationary solutions to the competitive Hele-Shaw problem are in a one-to-one correspondence with half-translation surfaces of Green's type.
\end{thm}
Theorem~\ref{thm:1-1} is proved in Section \ref{subsec:half}. 

\subsection{Interface erosion}\label{sect:introCIE}
\begin{figure}

\begin{center}
    
\begin{tikzpicture}

\filldraw[fill=blue!20!white, draw=black] (0,0) -- (0,0.5) -- (0.5,0.5) -- (0.5,1) -- (1,1) -- (1,0.5) -- (1.5,0.5) -- (2,0.5) -- (2.5,0.5) -- (2.5,0) -- (2.5,-0.5) -- (2,-0.5) -- (2,-1) -- (1.5,-1) -- (1.5,-1.5) -- (1,-1.5) -- (1,-1) -- (0.5,-1) -- (0.5,-0.5) -- (0,-0.5) -- (0,0);

\fill[blue!60!white] (1,0) -- (1.5,0) -- (1.5,-0.5) -- (1,-0.5) -- (1,0);

\filldraw[fill=green!20!white, draw=black] (2.5,0) -- (2.5,0.5) -- (2,0.5) -- (2,1) -- (3.5,1) -- (3.5,0.5) -- (4.5,0.5) -- (5,0.5) -- (5,-0.5) -- (4,-0.5) -- (4,-1) -- (4.5,-1) -- (4.5,-0.5) -- (5,-0.5) -- (5,-1.5)-- (3,-1.5) -- (3,-2) -- (2,-2) -- (2,-1.5) -- (1.5,-1.5) -- (1.5,-1) -- (2,-1) -- (2,-0.5) -- (2.5,-0.5) -- (2.5,0);

\fill[green!60!white] (3,-0.5) -- (3.5,-0.5) -- (3.5,-1) -- (3,-1) -- (3,-0.5);

\draw[green] (3.25,-0.75) -- (3.75,-0.75) -- (3.75,-0.25) -- (4.2,-0.25) -- (4.2,0.25) -- (4.3, 0.25) -- (4.3,-0.25) -- (4.75,-0.25) -- (4.75,-0.75);

\draw[thick,green,dashed] (4.5,-1) -- (5,-1);

\filldraw[fill=blue!20!white, draw=black] (6,0) -- (6,0.5) -- (6.5,0.5) -- (6.5,1) -- (7,1) -- (7,0.5) -- (7.5,0.5) -- (8,0.5) -- (8.5,0.5) -- (8.5,-1) -- (7.5,-1) -- (7.5,-1.5) -- (7,-1.5) -- (7,-1) -- (6.5,-1) -- (6.5,-0.5) -- (6,-0.5) -- (6,0);

\fill[blue!60!white] (7,0) -- (7.5,0) -- (7.5,-0.5) -- (7,-0.5) -- (7,0);

\filldraw[fill=green!20!white, draw=black] (8.5,0) -- (8.5,0.5) -- (8,0.5) -- (8,1) -- (9.5,1) -- (9.5,0.5) -- (10.5,0.5) -- (11,0.5) -- (11,-1) -- (10.5,-1) -- (10.5,-0.5) -- (10,-0.5) -- (10,-1) -- (10.5,-1) -- (11,-1) -- (11,-1.5)-- (9,-1.5) -- (9,-2) -- (8,-2) -- (8,-1.5) -- (7.5,-1.5) -- (7.5,-1) -- (8,-1) -- (8,-0.5) -- (8.5,-0.5) -- (8.5,0);

\fill[green!60!white] (9,-0.5) -- (9.5,-0.5) -- (9.5,-1) -- (9,-1) -- (9,-0.5);

\draw[blue] (7.25,-0.25) -- (7.25,0.25) -- (7.75,0.25) -- (7.75,-0.225) -- (8.25,-0.225) -- (8.25,-0.275) -- (7.75,-0.275) -- (7.75,-0.75) -- (8.25,-0.75);

\draw[thick,blue,dashed] (8,-1) -- (8.5,-1) -- (8.5,-0.5);

\end{tikzpicture}

\end{center}
\caption{Interface erosion.} 
\label{fig:6}
\end{figure}
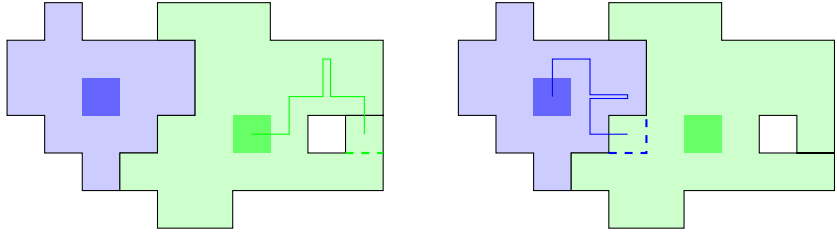

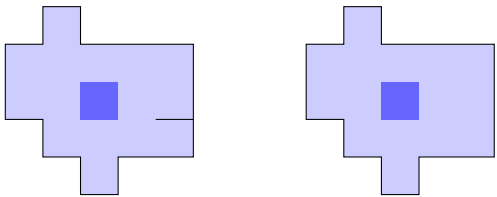
\begin{figure}

\begin{center}
    
\begin{tikzpicture}

\filldraw[fill=blue!20!white, draw=black] (0,0) -- (0,0.5) -- (0.5,0.5) -- (0.5,1) -- (1,1) -- (1,0.5) -- (1.5,0.5) -- (2,0.5) -- (2.5,0.5) -- (2.5,0) -- (2.5,-1) -- (1.5,-1) -- (1.5,-1.5) -- (1,-1.5) -- (1,-1) -- (0.5,-1) -- (0.5,-0.5) -- (0,-0.5) -- (0,0);

\fill[blue!60!white] (1,0) -- (1.5,0) -- (1.5,-0.5) -- (1,-0.5) -- (1,0);

\draw (2.5,-0.5) -- (2,-0.5);

 \filldraw[fill=blue!20!white, draw=black] (4,0) -- (4,0.5) -- (4.5,0.5) -- (4.5,1) -- (5,1) -- (5,0.5) -- (5.5,0.5) -- (6,0.5) -- (6.5,0.5) -- (6.5,-1) -- (5.5,-1) -- (5.5,-1.5) -- (5,-1.5) -- (5,-1) -- (4.5,-1) -- (4.5,-0.5) -- (4,-0.5) -- (4,0);

 \fill[blue!60!white] (5,0) -- (5.5,0) -- (5.5,-0.5) -- (5,-0.5) -- (5,0);

\end{tikzpicture}

\end{center}
\caption{Slits are removed.} 
\label{fig:7}
\end{figure}

We now introduce a random lattice version of the competitive Hele-Shaw flow. Consider the unit square lattice $\mathbb{L}$ viewed as a graph (i.e., the graph with vertices at $\mathbb{Z}^2$ and undirected edges connecting vertices at unit distance). Write $\mathbb{L}^*$ for the dual graph of $\mathbb{L}$ (i.e., the graph having the faces of $\mathbb{L}$ as its vertices and edges connecting vertices at unit distance).
\begin{df}[Interface erosion] \label{df:IE}
Let $(D^i_0,\gamma^i_0,d^i)$ be an admissible droplet configuration on $\CC$ such that all droplet interfaces lie in $\mathbb{L}$. For each $i$ and  each source point $z^i_j$ in the source divisor $d^i=\sum_j a^i_j z^i_j$ we identify $z^i_j$ with the square in which it lies (if it lies in several squares we choose one). We thus think $z^i_j$ of as a vertex in $\mathbb{L}^*$, and we associate to it an independent Poisson clock with rate $a^i_j$. Suppose the clock at $z^i_j$ rings at time $t$. The interfaces at time $t$ are then obtained from those at time $t-$ as follows. We start a random walk on $\mathbb{L}^*$ from $z^i_j$. We stop the walk when it crosses the droplet interface $\gamma^i_{t-}$ at an edge $\mathbf{e}$. If the end square is a source square all interfaces are left unchanged. Otherwise, all droplet interfaces $\gamma^j_{t-}$ that pass through $\mathbf{e}$ are redirected so that they instead go around the end square, see Figure \ref{fig:6}. If an interface develops a slit, then that part of the interface is removed, see Figure \ref{fig:7}. We discuss below how to define the model on more general surfaces.

\begin{rem}
    Instead of using Poisson clocks one could just as well let the random walk beginning at $z^i_j$ start at regular intervals of length $1/a^i_j$, and if then two or more random walks are to start at the same time, we perform them in say lexicographic order of the indices $(i,j)$.
\end{rem}

\end{df}

Interface erosion is a slight variation of Propp's competitive erosion model \cite{GLPP}, see Section~\ref{sect:discretemodels}. In both models internal DLA-like clusters compete to grow but the essential difference is that the topology of clusters is preserved in interface erosion whereas it is not in competitive erosion. However, in suitable circumstances we expect long term-limits to be the same. (Roughly speaking, this should be the case if interface erosion is started with the ``correct'' topology, i.e., the one which competitive erosion eventually approaches.) We think of interface erosion as a regularized version of the competitive Hele-Shaw flow, analogous to how internal DLA and DLA relates to classical Hele-Shaw with injection and suction, respectively. We expect that long-term limits of interface erosion, appropriately rescaled, are generically described by stationary solutions to the competitive Hele-Shaw flow, see Conjecture~\ref{conj:comper} in Section~\ref{sect:discretemodels}. In fact, describing possible scaling limits of competitive erosion in some generality was part of the motivation for the present paper. 

\subsubsection*{Simulations of interface erosion}

\begin{figure}[t]

\centering

\begin{subfigure}{.25\textwidth}

    \includegraphics[width=.8\linewidth]{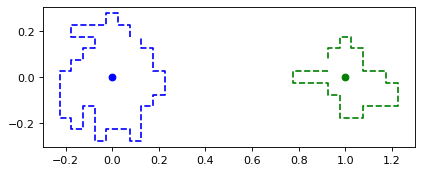}
    \caption{t=0.1}
    \label{fig:0.1}
\end{subfigure}%
\begin{subfigure}{.25\textwidth}
    \includegraphics[width=.8\linewidth]{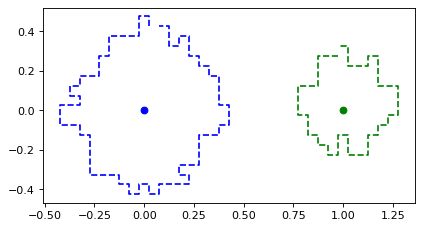}
    \caption{t=0.25}
    \label{fig:0.25}
\end{subfigure}%
\begin{subfigure}{.25\textwidth}
    \includegraphics[width=.8\linewidth]{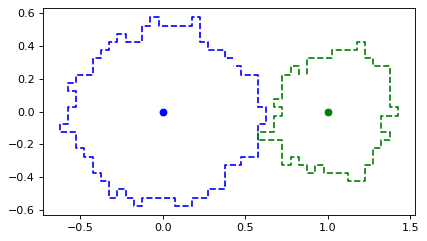}
    \caption{t=0.5}
    \label{fig:0.5}
\end{subfigure}%
\begin{subfigure}{.25\textwidth}
    \includegraphics[width=.8\linewidth]{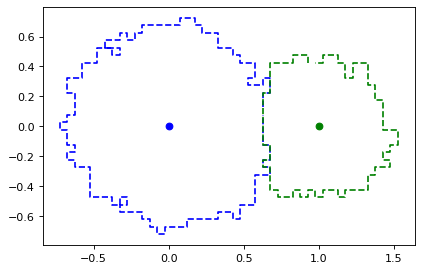}
    \caption{t=0.75}
    \label{fig:0.75}
\end{subfigure}

\begin{subfigure}{.25\textwidth}
    \includegraphics[width=.8\linewidth]{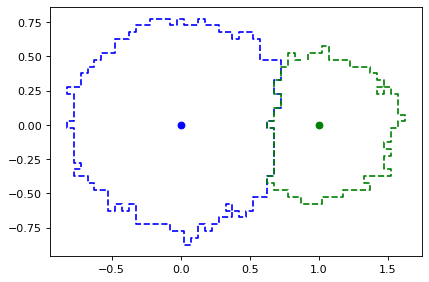}
    \caption{t=1}
    \label{fig:1.0}
\end{subfigure}%
\begin{subfigure}{.25\textwidth}
    \includegraphics[width=.8\linewidth]{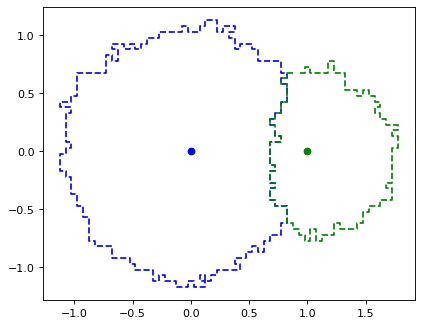}
    \caption{t=2}
    \label{fig:2.0}
\end{subfigure}%
\begin{subfigure}{.25\textwidth}
    \includegraphics[width=.8\linewidth]{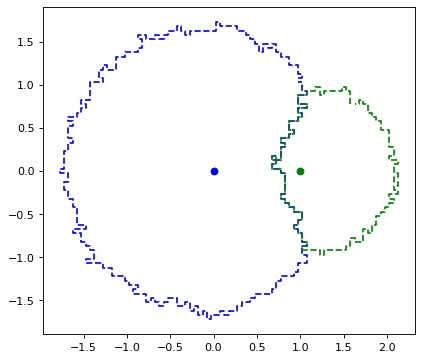}
    \caption{t=5}
    \label{fig:5.0}
\end{subfigure}%
\begin{subfigure}{.25\textwidth}
    \includegraphics[width=.8\linewidth]{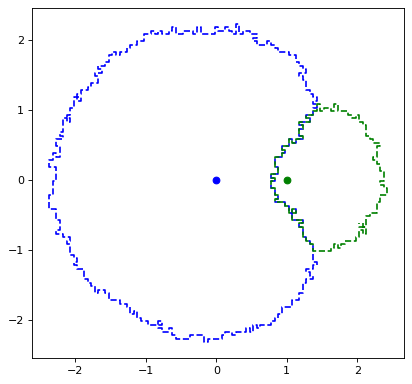}
    \caption{t=10}
    \label{fig:10.0}
\end{subfigure}

\begin{subfigure}{.25\textwidth}
    \includegraphics[width=.8\linewidth]{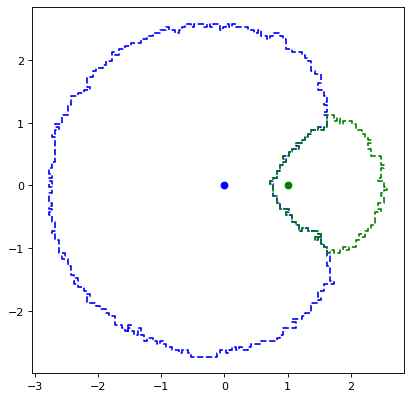}
    \caption{t=15}
    \label{fig:15.0}
\end{subfigure}%
\begin{subfigure}{.25\textwidth}
    \includegraphics[width=.8\linewidth]{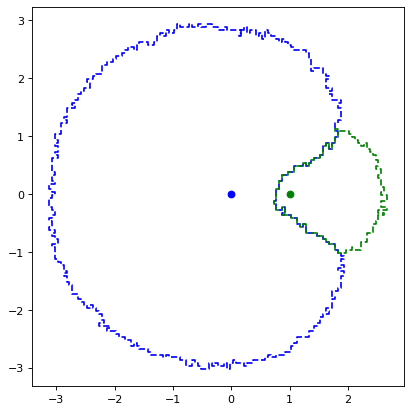}
    \caption{t=20}
    \label{fig:8000}
\end{subfigure}%
\begin{subfigure}{.25\textwidth}
    \includegraphics[width=.8\linewidth]{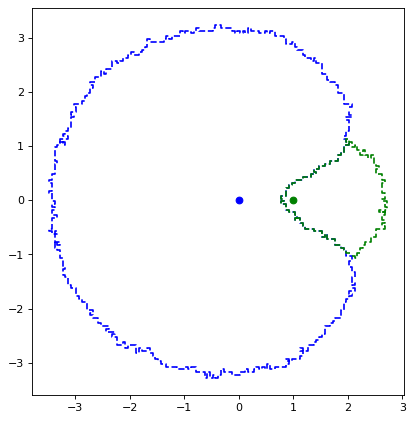}
    \caption{t=25}
    \label{fig:25.0}
\end{subfigure}%
\begin{subfigure}{.25\textwidth}
    \includegraphics[width=.8\linewidth]{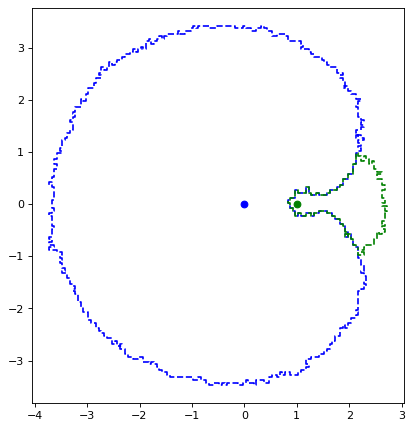}
    \caption{t=30}
    \label{fig:30.0}
\end{subfigure}

\begin{subfigure}{.25\textwidth}
    \includegraphics[width=.8\linewidth]{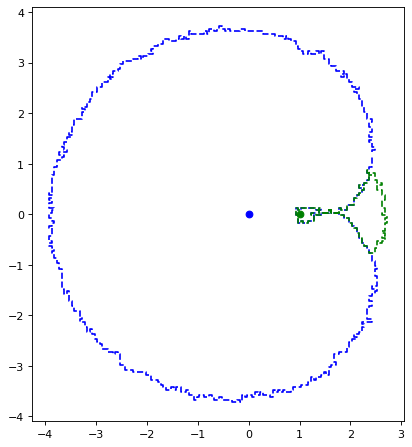}
    \caption{t=35}
    \label{fig:35.0}
\end{subfigure}%
\begin{subfigure}{.25\textwidth}
    \includegraphics[width=.8\linewidth]{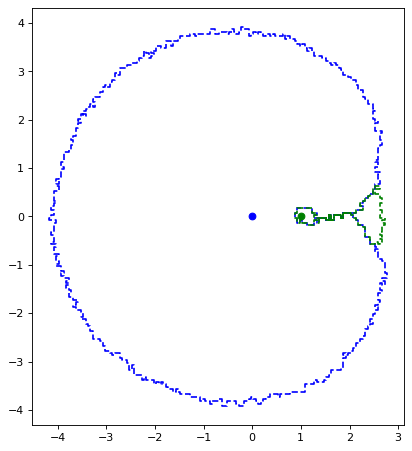}
    \caption{t=40}
    \label{fig:40.0}
\end{subfigure}%
\begin{subfigure}{.25\textwidth}
    \includegraphics[width=.8\linewidth]{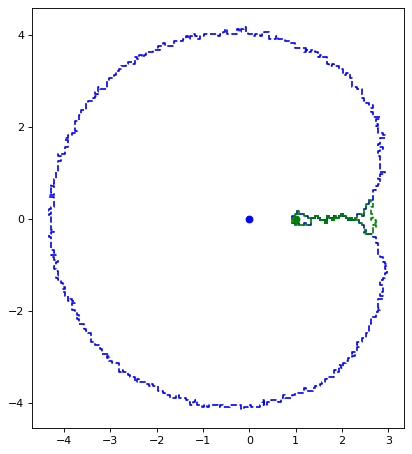}
    \caption{t=45}
    \label{fig:45.0}
\end{subfigure}%
\begin{subfigure}{.25\textwidth}
    \includegraphics[width=.8\linewidth]{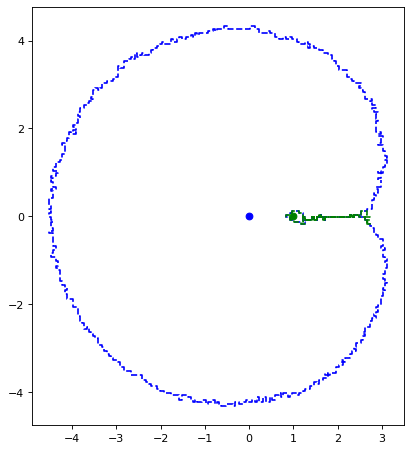}
    \caption{t=50}
    \label{fig:50.0}
\end{subfigure}

\caption{Simulation of interface erosion on $\CC$ with source divisors $d^1=2\cdot 0$ and $d^2=1\cdot 1$ (starting from small circles around the source points) using a square grid of mesh size $1/20$. It seems plausible that the discrete droplets approach the trajectories given by the quadratic differential \eqref{example:teichmuller} which is the one constructed in Proposition~\ref{thm2} for these initial data.}

\label{fig:sim}

\end{figure}

\begin{figure}

\centering

\begin{subfigure}{.5\textwidth}
\centering
    \includegraphics[width=0.67\linewidth]{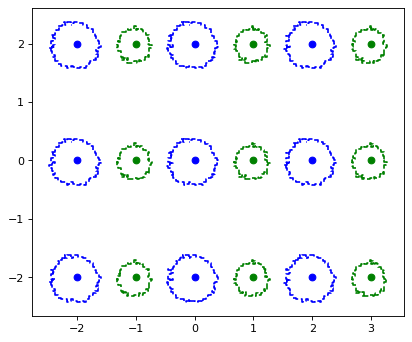}
    \caption{t=0.25}
    \label{fig:t0.25}
\end{subfigure}%
\begin{subfigure}{.5\textwidth}
\centering
    \includegraphics[width=0.67\linewidth]{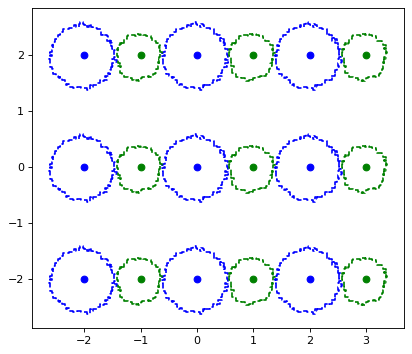}
    \caption{t=0.5}
    \label{fig:t0.5}
\end{subfigure}

\begin{subfigure}[t]{.5\textwidth}
\centering
    \includegraphics[width=0.67\linewidth]{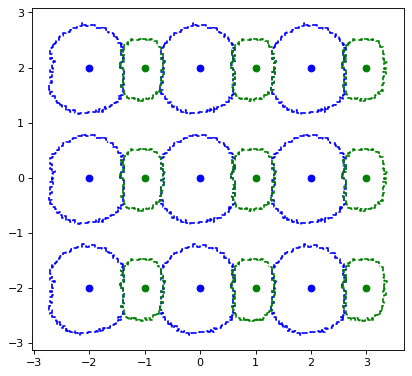}
    \caption{t=1}
    \label{fig:t1}
 \end{subfigure}%
 \begin{subfigure}[t]{.5\textwidth}
 \centering
    \includegraphics[width=0.67\linewidth]{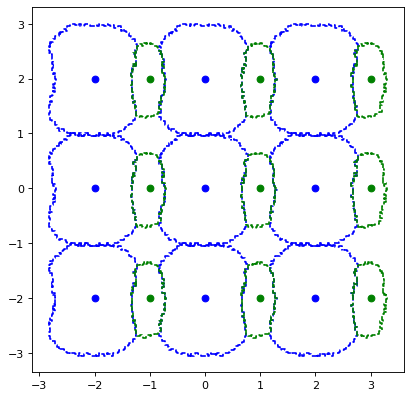}
    \caption{t=2}
    \label{fig:t2}
\end{subfigure}

\begin{subfigure}{.5\textwidth}
\centering
    \includegraphics[width=0.67\linewidth]{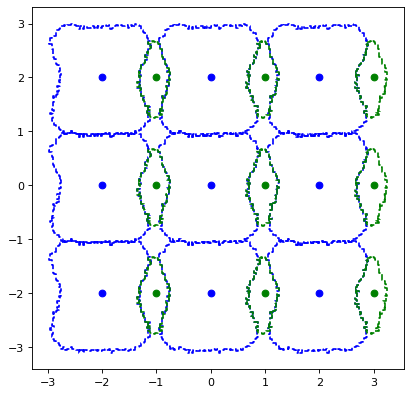}
    \caption{t=3}
    \label{fig:t3}
\end{subfigure}%
\begin{subfigure}{.5\textwidth}
\centering
    \includegraphics[width=0.67\linewidth]{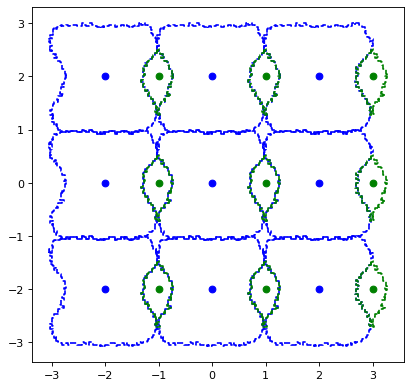}
    \caption{t=4}
    \label{fig:t4}
\end{subfigure}

\begin{subfigure}{.5\textwidth}
\centering
    \includegraphics[width=0.67\linewidth]{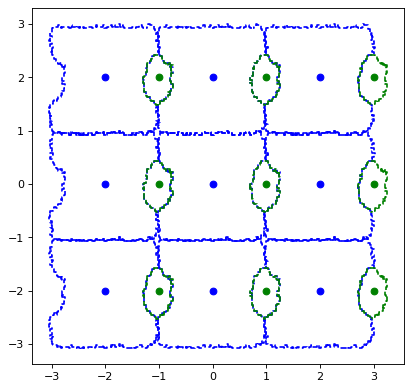}
    \caption{t=5}
    \label{fig:t5}
\end{subfigure}%
\begin{subfigure}{.5\textwidth}
\centering
    \includegraphics[width=0.677\linewidth]{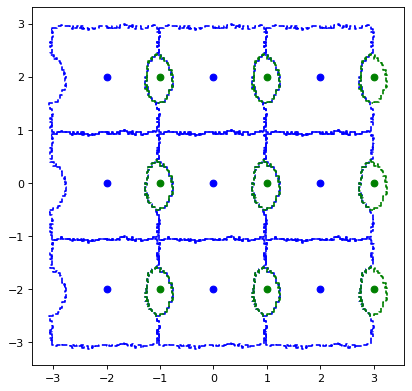}
    \caption{t=6}
    \label{fig:t6}
\end{subfigure}

\caption{Simulation of interface erosion on the torus with fundamental domain $[-0.5,1.5]\times [-1,1]$ with the induced flat metric and with sources at $(0,0)$ and $(1,0)$ of strength $2$ and $1$, respectively. Mesh size $1/40$. For large $t$ interfaces are very close to the trajectories of the quadratic differential constructed in Proposition~\ref{thm2} with the given initial data.}
\label{fig:sim2b}

\end{figure}

Simulations of the interface erosion model lend some support for Conjecture \ref{conj:cHS} and Conjecture \ref{conj:comper}. Admitting the latter we can think of interface erosion as a way to simulate solutions to the competitive Hele-Shaw problem, in particular in the long-term limit. Figure~\ref{fig:sim} and Figure~\ref{fig:sim2b} show the result of two such simulations. 

\subsection{Further remarks}

\paragraph{Competitive Hele-Shaw flow as a gradient flow.} Consider two droplets $D,D^*$ with divisors $d,d^*$ in $\hat{\mathbb{C}}$ separated by a smooth bounded Jordan curve $\gamma$. We write $I(\gamma, d,d^*)$ for the reduced Green's energy. (See Section~\ref{sect:jordan}.) In the case when $d=1 \cdot 0, d^* = 1 \cdot \infty$, we have $I(\gamma, d,d^*) =\log r_D(0) + \log r_{D^*}(\infty)$. The right-hand side has an interpretation as a K\"ahler potential for a canonical metric on the universal Teichm\"uller curve (roughly speaking, a space of normalized quasicircles in $\hat{\mathbb{C}}$ with a marked point in the Riemann sphere), see \cite{TT}. This role is played by the universal Liouville action/Loewner energy in Weil-Petersson Teichm\"uller space (a space of normalized quasicircles), see \cite{W,TT} for background and definitions. In \cite{BMVW} a gradient flow in universal Teichm\"uller space with respect to the Loewner energy is considered.  Proposition~\ref{prop:gradient} shows that the competitive Hele-Shaw flow (in the sense of Hadamard variation) decreases $-I(\gamma, d,d^*)$, and with a suitable interpretation (e.g., as in Chapter~6.3 of \cite{walker}) one can view the competitive Hele-Shaw flow as a gradient flow for the reduced Green's energy. But it is not yet clear to us whether this also holds in a stronger sense as in \cite{BMVW}. Let us also mention that by adding a small multiple of the Loewner energy one formally obtains a natural notion of viscosity solution. We will study these and related questions elsewhere \cite{duse}.

\begin{figure}

\begin{center}

\begin{tikzpicture}

\fill[blue!20!white] (-1,1.5) -- (0,0) -- (-1,-1.5) -- (-1,1.5);

\fill[green!20!white] (-1,1.5) -- (0,0) -- (1,0) -- (2,1.5) -- (-1,1.5);

\fill[red!20!white] (-1,-1.5) -- (0,0) -- (1,0) -- (2,-1.5) -- (-1,-1.5);

\fill[yellow!20!white] (2,1.5) -- (2,-1.5) -- (1,0) -- (2,1.5);

\fill[blue!20!white] (4,1.5) -- (5.5,0.5) -- (5.5,-0.5) -- (4,-1.5) -- (4,1.5);

\fill[green!20!white] (4,1.5) -- (5.5,0.5) -- (7,1.5) -- (4,1.5);

\fill[red!20!white] (4,-1.5) -- (5.5,-0.5) -- (7,-1.5) -- (4,-1.5);

\fill[yellow!20!white] (7,1.5) -- (5.5,0.5) -- (5.5,-0.5) -- (7,-1.5) -- (7,1.5);

\draw  (-1,1.5) -- (0,0) -- (-1,-1.5);

\draw (-1,1.5) -- (0,0) -- (1,0) -- (2,1.5);

\draw (-1,-1.5) -- (0,0) -- (1,0) -- (2,-1.5);

\draw  (2,1.5) -- (1,0) -- (2,-1.5);

\draw (4,1.5) -- (5.5,0.5) -- (5.5,-0.5) -- (4,-1.5);

\draw (4,1.5) -- (5.5,0.5) -- (7,1.5);

\draw (4,-1.5) -- (5.5,-0.5) -- (7,-1.5);

\draw (7,1.5) -- (5.5,0.5) -- (5.5,-0.5) -- (7,-1.5);

\path (-0.6,0) node {$D^1$};

\path (0.6,0.75) node {$D^2$};

\path (1.6,0) node {$D^3$};

\path (0.6,-0.75) node {$D^4$};

\path (4.75,0) node {$D^1$};

\path (5.6,1.1) node {$D^2$};

\path (6.25,0) node {$D^3$};

\path (5.6,-1.1) node {$D^4$};

\end{tikzpicture}

\end{center}
    
\caption{A transition between droplet configurations corresponding to a wall-crossing.} 
\label{fig:wc}
\end{figure}
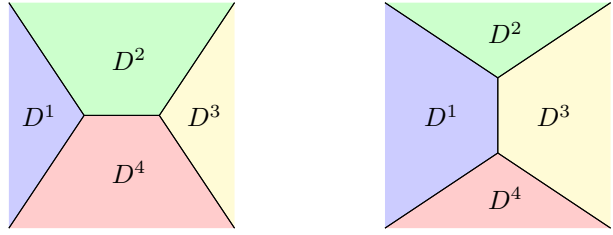

\paragraph{Bridgeland stability conditions and wall-crossing.} It is interesting to note that the kind of quadratic differentials that appear in this paper (i.e., meromorphic quadratic differentials with second order poles) also play a role in the theory of Bridgeland stability conditions, as shown by Bridgeland-Smith \cite{BS}.  Wall-crossing, a phenomenon of great importance in that context, occurs when four droplets that initially are in a configuration as on the left in Figure \ref{fig:wc} transition to a configuration as on the right. That is, while first $D^2$ and $D^4$ share a common boundary, after the transition it is $D^1$ and $D^3$ that share a common boundary. Note that this often can be achieved e.g. by increasing the source strengths of $d^1$ and $d^3$ while decreasing the source strengths of $d^2$ and $d^4$. Alternatively one can move the source points of $d^1$ and $d^3$ closer together.

\subsection*{Acknowledgements}
F.V. acknowledges support from the Knut and Alice Wallenberg Foundation, the Göran Gustafsson Foundation for Research in Natural Sciences and Medicine, and the Swedish Research Council. D.W.N.\ has been supported by the Göran Gustafsson Foundation for Research in Natural Sciences and Medicine, and the Swedish Research Council. We thank Bj\"orn Gustafsson and Yilin Wang for discussions and useful comments on an earlier version the paper and Erik Duse and Julius Ross for discussions. We are also grateful to Steffen Rohde for discussions and for constructing a stationary Jordan curve which led us to the link to lemniscates.    
\section{Preliminaries}

\subsection{Basic definitions} \label{subsec:prel}
Let $\Sigma$ be a Riemann surface. A (positive) \emph{divisor} $d = \sum_{j=1}^n a_j z_j$ on $\Sigma$ is a finite linear combination of points $z_j \in \Sigma$ with (positive) real coefficients $a_j$. The weight of $d$ is defined by $|d| = \sum_{j=1}^n a_j$ and the support of $d$ is defined by $\textrm{supp}(d) = \{z_j: a_j\neq 0\}$. We say that two divisors are disjoint if their supports are disjoint. The sum of two (disjoint) divisors is defined in the obvious way by $\sum_{j=1}^n a_j z_j+\sum_{j=1}^m b_j w_j = \sum_{j=1}^{m+n}c_j v_j$, where $c_j=a_j, \, v_j=z_j,  j=1,\ldots, n$ and $c_j=b_{j-n},\, v_j=w_{j-n},\, j=n+1, \ldots, m$.
If $f:\Sigma \to \Sigma$ is any function, we set $f(d) = \sum_{j} a_j f(z_j)$. 
  
We will often consider $n$-tuples of divisors $d=(d^i)_{i=1}^n$ and the definitions of weight and support extend naturally. 

Let $D$ be a planar finitely connected Jordan domain (say). The conformally invariant Green's function for $D$ with pole at $w$ is defined by
\[
G_{D}(z,w) = G_{D,w}(z) = \log |z-w|^{-1} + h_{D,w}(z) ,
\]
where $z \mapsto h_{D,w}(z)$ is the solution to the Dirichlet problem with boundary data $\zeta \mapsto \log|\zeta-w|$. Note that $\log |z-w|^{-2}  =\log(z-w)^{-1} + \log(\bar z-\bar w)^{-1}$, so
\[
\partial_z G_{D,w}(z) = \frac{-1/2}{z-w} + \partial_z h_{D,w}(z).
\]

In the unit disc we have
\[
G_{\mathbb{D},w}(z) = \log \left|\frac{1-z \bar{w}}{z-w}\right|, \quad \partial_zG_{\mathbb{D},w}(z) = -\frac{1}{2}\left(\frac{1}{z-w}+ \frac{\bar w}{1-z \bar w}\right).
\]
The Green's function is also defined on a hyperbolic Riemann surface $\Sigma$. If the surface is simply connected, there is a conformal bijection $f:\Sigma \to \mathbb{D}$ and 
\[
G_\Sigma(x,y) = G_{\mathbb{D}}(f(x),f(y)). 
\]
 If the boundary of $D$ is sufficiently smooth, then harmonic measure and arc length are absolutely continuous and the Poisson kernel is defined by
 \[
 P_D(z,\zeta) = - \frac{\partial}{\partial{n_\zeta}} G_D(z,\zeta) > 0,
 \]
 where $n_\zeta$ is the outward pointing normal. (Note that we do not include the factor $1/2\pi$.) In particular,
 \[
 P_{\mathbb{D}}(z,\zeta)= \frac{1-|z|^2}{|z-\zeta|^2}.
 \]

Let $D \subset \Sigma$ be a simply connected open set. Then the Green's function for $D$ exists if and only if there is a conformal bijection $f: D \to \mathbb{D}$ and in this case for all $x,y \in D$, $G_D(x,y) = G_{\mathbb{D}}(f(x), f(y))$. $D$ is a Jordan domain in $\Sigma$ if it is connected and simply connected, $\partial D$ is a Jordan curve, and $\overline{D} \subset \Sigma$.

If $\gamma$ is a Jordan curve on $\Sigma$ homotopic to a point, then $\gamma$ is the boundary of a topological disc $D$ which is conformally equivalent to $\mathbb{D}$ and in particular, the Green's function for $D$ exists. Unless $\Sigma = \hat{\m C}$, $D$ is uniquely determined. Similarly, if $\gamma_1, \gamma_2$ are two disjoint and freely homotopic Jordan curves on $\Sigma$, then the curves form the two boundary components of a topological annulus which is conformally equivalent to a round annulus in $\mathbb{C}$. See \cite{strebel-book} Section~2.4 and 2.5.

Given a divisor $d=\sum_{k=1}^n a_k z_k$ in $D$, we write
\[
G_{D,d}(z) = \sum_{k=1}^n a_k G_{D,z_k}(z).
\]

\subsubsection*{A compactification of a punctured Riemann surface}

Let $\Sigma=\tilde{\Sigma}\setminus\{p_1,...,p_m\}$ where $\tilde{\Sigma}$ is a closed Riemann surface. For each puncture $p_i$ let $f_i:\mathbb{D}\to \tilde{\Sigma}$ be a conformal injection such that $f_i(0)=p_i$, and let $D_i:=f_i(\mathbb{D})$ and $D_i^*:=f_i(\mathbb{D}^*)$ (here $\mathbb{D}^*$ denotes the unit disk punctured at $0$). Note that the map $$g(z):=\frac{(1+|z|)z}{2|z|}$$ between $\mathbb{D}^*$ and the annulus $A:=\{z: 1/2<|z|<1\}$ is a diffeomorphism. Let $\bar{A}_i:=\{z: 1/2\leq|z|<1\}$. We now define $\mathring{\Sigma}$ to be the compact differentiable surface with boundary we get by first taking the disjoint union $\Sigma \cup_{i=1}^m \bar{A}_i$ and then for each $i$ identifying $D_i^*$ with $A_i:=\{z: 1/2<|z|<1\}\subseteq \bar{A}_i$ via the diffeomorphism $g\circ f_i^{-1}$. Thus there is a boundary circle $\{z: |z|=1/2\}\subseteq \bar{A}_i$ in $\mathring{\Sigma}$ for each puncture $p_i$.

The reason for letting the droplet interfaces lie in $\mathring{\Sigma}$ (rather than in $\Sigma$ or $\tilde{\Sigma}$) is that we want to allow a droplet interface $\gamma$ to ``touch'' a puncture but at the same time it should not be possible to continuously deform $\gamma$ across a puncture. Note however that there is a natural continuous projection map from $\mathring{\Sigma}$ to $\tilde{\Sigma}$, mapping each boundary circle to its associated puncture. Thus a curve in $\mathring{\Sigma}$ can be projected to a curve in $\tilde{\Sigma}$, and we can thus think of a droplet interface $\gamma$ as a curve in $\tilde{\Sigma}$ together with some additional information of how it goes round any puncture it crosses. 

We will also use a related surface $\Sigma_+$ defined as follows. Let $\mathbb{D}^*_i:=\{z: 0<|z|<1\}$. We now define $\Sigma_+$ to be the differentiable surface (without boundary) we get by first taking the disjoint union $\Sigma \cup_{i=1}^m \mathbb{D}^*_i$ and then for each $i$ identifying $D_i^*$ with $A_i:=\{z: 1/2<|z|<1\}\subseteq \mathbb{D}^*_i$ via the diffeomorphism $g\circ f_i^{-1}$. Clearly $\mathring{\Sigma}\subseteq \Sigma_+$.

\subsubsection*{Quadratic differentials and vertical foliations}

A meromorphic \emph{abelian} differential $\psi$ on $\Sigma$ is a collection of local meromorphic functions $\psi_\alpha$ obeying the following transformation law under coordinate change
\[\psi_\alpha(z_\alpha) dz_\alpha = \psi_\beta(z_\beta) dz_\beta,
\] while similarly a meromorphic \emph{quadratic} differential $\vp$ on $\Sigma$ is a collection of meromorphic functions $\vp_\alpha$ that follow the transformation law
\[
\vp_\alpha(z_\alpha) dz_\alpha^2 = \vp_\beta(z_\beta) dz_\beta^2.
\]

The zeroes and poles of $\vp$ are called critical points.

A quadratic differential $\vp$ on $\Sigma$ gives rise to a so-called vertical foliation of $\Sigma$. Let $\gamma = \gamma(t), \, t \in I,$ be a differentiable curve on $\Sigma$ and consider $t \mapsto z \circ \gamma(t) = :z(t)$ for a local coordinate $z$. Suppose that $\gamma$ does not pass through a critical point of $\vp$. If
$$\arg( \vp(z(t)) z'(t)^2) = 2\theta, \quad \theta\in[0,\pi), \quad t \in I,$$
we say that $\gamma$ is a flow line of angle $\theta$ for $\vp$. If $\theta = 0$ or $\theta =\pi/2$, $\gamma$ is said to be a horizontal flow line and vertical flow line, respectively. Another way to see it is that if $\sqrt{\vp}$ is a local square root and $z_0$ is a local reference point then the vertical foliation is given by the level sets of $\int_{z_0}^z \textrm{Re}(\sqrt{\vp})$. A \emph{trajectory} is a maximal flow line. Note that a flow line of angle $\theta$ for $\vp$ is a horizontal flow line for $e^{-i 2\theta } \vp$.  

\emph{Example.} Let $d=\sum_{j=1}^n a_j z_j$ be a divisor in $\mathbb{D}$ and consider the quadratic differential 
\[
\varphi(z)dz^2 = (\partial_z G_{\mathbb{D},d}(z))^2dz^2.
\]
Suppose $t\mapsto z(t)$ parametrizes a smooth part of an equipotential (level line) for $z \mapsto \psi(z):=\sum_{j=1}^n a_j G_{\mathbb{D}}(z, z_j)$. Then 
\[
0=\frac{d}{dt}\psi(z(t)) = 2 \textrm{Re} \sum_{j=1}^n a_j \partial_z G_{\mathbb{D}}(z(t), z_j)z'(t).
\]
Hence, where $z(t)$ is differentiable,
\[
 \arg \left( \varphi(z(t))z'(t)^2 \right) \equiv \pi
\]
and we see that $t \mapsto z(t)$ determines a vertical trajectory for the quadratic differential $\varphi$. Note that $\psi = 0$ on $\partial \mathbb{D}$ and since the weights $a_j$ are positive, $\psi$ is superharmonic and $ > 0$ in $\mathbb{D}$. Therefore there is some maximal $\lambda_0 > 0$ (depending only on $d$) such that the $\lambda$-level lines of $\psi$, $\lambda \in (0,\lambda_0)$, are Jordan curves which each separates $\partial \mathbb{D}$ from $\textrm{supp}(d)$. As $\lambda$ is varied, these curves sweep out an ``outermost'' ring domain in $\mathbb{D}$ which clearly is also a characteristic ring domain for the vertical trajectories of $\varphi$. As $\lambda$ is further increased, the level lines sweep out further characteristic ring domains and, eventually, as $\lambda \to \infty$, once-punctured topological disks near the individual points in $\textrm{supp}(d)$.

\subsubsection*{Translation and half-translation surfaces}
A (half-) translation surface is a collection of polygons $P_i$ in $\CC$ with sides pairwise identified via translations (and/or minus-translations $z \mapsto -z+b$). Note that if $P_i$ has a side $s_i$ which is identified with a side $s_j$ in $P_j$ via a translation $z\mapsto z+a$, we require that $P_i+a$ and $P_j$ lie on different sides of $s_j=s_i+a$. Similarly, if $P_k$ has a side $s_k$ which is identified with a side $s_l$ in $P_l$ via a minus-translation $z\mapsto -z+b$, $-P_k+b$ and $P_l$ should lie on different sides of $s_l=-s_k+b$.

The most basic example of a translation surface is a single rectangle, say $S=\{z: a\leq \textrm{Re}(z)\leq b, c\leq \textrm{Im}(z)\leq d,\}$, with the left and right side identified via the translation $z\mapsto z+b-a$ and the top and bottom side identified via the translation $z\mapsto z+(c-d)i$.

A translation surface gives rise to a Riemann surface together with an abelian differential, given by $dz$ on any given polygon. In the example above we of course get a torus, with its nonvanishing abelian differential. Going in the other direction, any Riemann surface with an abelian differential can be represented by a translation surface, and given two such representations one can go from one to the other via a simple cut and gluing operation.

A simple example of a half-translation is given by the rectangle $R=\{z: 0\leq \textrm{Re}(z)\leq 1, -1\leq \textrm{Im}(z)\leq 1,\}$ where we identify the top and bottom side via the translation $z\mapsto z+2i$, but in contrast to the earlier example we now split the left side $[-i,i]$ into two equal parts: $[-i,0]$ and $[0,i]$ and identify them via the minus-translation $z\mapsto -z$. Similarly we split the right side $[1-i,1+i]$ into two equal parts: $[1-i,1]$ and $[1,1+i]$ and identify them via the minus-translation $z\mapsto -z+2$.

A half-translation surface gives rise to a Riemann surface together with a quadratic differential, given by $dz^2$ on any given polygon. In the example above we get $\hat{\CC}$, and the quadratic differential on $\hat{\CC}$ that corresponds to $dz^2$ on $R$ will have simple poles at the four points in $\hat{\CC}$ that corresponds to the four points $0,i,1,1+i\in R$ (note that $-i$ and $1-i$ corresponds to the same points as $i$ and $1+i$). Note also that the flat metric on $R$ induces a metric on $\hat{\CC}$ which is flat except at these four points where it will have conical singularities, all with cone angle $\pi$.

Going in the other direction, any Riemann surface with a quadratic differential can be represented by a half-translation surface, and given two such representations one can go from one to the other via a simple cut and gluing operation. On the associated half-translation surface the vertical foliation is simply given by the level sets of $\textrm{Re}(z)$. 

For an introduction to the theory of translation and half-translation surfaces see \cite{Zor}.

Sometimes one has a collection of polygons where some but not all sides are pairwise identified via translations (and/or minus-translations). We will in this paper call such a surface a \emph{partial} (half-) translation surface.

\subsection{Reduced Green's energy} 

Let $\Sigma$ be a Riemann surface of finite type. Let $D$ be a simply connected domain in $\Sigma$ and $d = \sum_{j=1}^n a_j z_j$ a divisor with support in $D$. Given $z_j$, choose a local coordinate and for small $r$ let $B_r$ be the preimage in $\Sigma$ of the disc of radius $r$ around $0$. Let $D_j$ be the component of $D$ containing $z_j$ and let $A_r = D_j\smallsetminus B_r$. This is a topological annulus for all sufficiently small $r > 0$ and we define the reduced modulus of $D\smallsetminus\{z_j\}$ by
\[M_D(z_j)=\lim_{r \to 0} 2\pi \left(\textrm{Mod}(A_r) - \frac{1}{2\pi} \log r^{-1}.\right)\]
Here $\textrm{Mod}(A_r)$ is the conformal modulus of the annulus $A_r$. Note that, as opposed to the conformal modulus, the reduced modulus is not conformally invariant  and depends on the choice of coordinate.
In the plane we have equivalently
\[
M_D(z_j) = \lim_{z \to z_j}(G_{D, z_j}(z) -\log|z-z_j|^{-1}).
\]

We define the \emph{reduced Green's energy} of $(D,d)$ by
\begin{align}\label{def:Green-Energy}
I(D, d) = \sum_{1 \le j \neq k \le n} a_ja_kG_D(z_j,z_k) + \sum_{j=1}^n a_j^2M_{D}(z_j).
\end{align}

If $D$ is a droplet with possibly several components $E^j$, we define 
\[I(D,d) = \sum_{j=1}^mI(E^j, d^j)\] where $d^j = d \cap E^j$.
 
 Since all terms in \eqref{def:Green-Energy} have this property, we see that $I(D, d)$ is monotone increasing in $D$ for $d$ fixed. In the plane, we have the following interpretation. The off-diagonal terms in \eqref{def:Green-Energy} correspond to the electrostatic potential energy of the ensemble of charges $\{a_j\}$ at the marked points $\{z_j\}$ in $D$ with the boundary $\partial D$ grounded.  We usually think of the energy \eqref{def:Green-Energy} as being attached to the domain (which we will vary) and not to the marked points (which always stay fixed). Note that only the harmonic part of $I(D,d)$ changes with the domain.

It is clear that if $f : D \to f(D)$ is a conformal bijection, then
\[ I(f(D), f(d)) = I(D, d) +  \sum_{j=1}^n a_j^2\log|\partial_z (w\circ f)(z_j)|.\]

Let $\mathcal{D} = (D^i,\gamma^i,d^i)_{i=1}^n$ be an admissible droplet configuration on $\Sigma$, and assume that we have chosen local coordinates around all source points. Then we define the \emph{reduced Green's energy} of the droplet configuration as 
\begin{align}\label{def:green-energy-droplet-config}
I(\mathcal{D}) = \sum_{i=1}^n I(D^i, d^i).
\end{align}
It is immediate from the definition that $I(\mathcal{D}) > - \infty$. 
 
\subsection{Hadamard- and quasiconformal variation}
Let $\gamma$ be a $C^2$ Jordan curve in $\hat{\mathbb{C}}$ bounding the domain $D$ and let $\nu(z)$ be a $C^1$ function of arc-length along $\gamma$. (We do not consider the most general setting with respect to regularity here.) Let $n(z)$ be the normal vector of $\gamma$ in the outward direction with respect to $D$. Hadamard's classical variational formula describes the first variation of the Green's function of $D$ under the variation of $\gamma$ by the vector field $\nu(z)n(z)$. For $t$ small enough $s \mapsto \gamma_t(s) = \gamma(s) + t\nu(\gamma(s))n(\gamma(s))$ defines a Jordan curve with inner domain $D_t$ also containing the points $z,w$. Note that we do not assume $\nu$ is positive. See, e.g., Appendix~3 of \cite{courant} for this version of Hadamard's formula.
 \begin{lemma}\label{lemma:hadamard}
 For all $t$ in a sufficiently small neighborhood of $0$, 
 \begin{align}\label{eq:hadamard}
  G_{D_t,z_1}(z)-G_{D,z_1}(z) =  \frac{t}{2\pi}\int_\gamma P_{D,z_1}(\zeta)P_{D,z}(\zeta)\nu(\zeta) |d\zeta| + O(|t|^2),
 \end{align}
 where the error term is uniform for $z,z_1$ in a given compact subset of $D$.
    \end{lemma}
    We may write \eqref{eq:hadamard} concisely as \[\delta_\nu G_D(z,w) =  \frac{1}{2\pi}\int_\gamma P_{D,z}(\zeta)P_{D,w}(\zeta)\nu(\zeta) |d\zeta|\]
    and we will use this notation below for Hadamard variations.
Let $h_{D}(z,w)=G_D(z,w) - \log|z-w|^{-1}$ be the harmonic part of the Green's function. Then $\delta_\nu h_{D}(z,w) =  \delta_\nu G_{D}(z,w)$ since the varied functions differ by a function independent of $\gamma$. This implies that the Hadamard variation of the reduced modulus is given by
\[
\delta_\nu M_D(z) = \frac{1}{2\pi} \int_{\partial D}P_{D,z}(\zeta)^2 \nu(\zeta) |d\zeta|.
\]
The meaning of this formula is as in Lemma~\ref{lemma:hadamard}.

We shall also make use of quasiconformal variation of the Green's function. While it is possible to view the Hadamard variation as a special case of the quasiconformal variation below, we choose to present the former in the classical way. The following lemma follows from \cite{sontag}. Here we do not make any assumptions on (the regularity of) $D$ besides being simply connected.
\begin{lemma}\label{lem:variation}
Let $f_t$ be a quasiconformal map of $\mathbb{C}$ with dilatation $\mu_t$ with compact support. Suppose $\|\mu_t\|_\infty = O(|t|)$ and $F_t(z) =z + O(|t|)$ as $t \to 0+$. Let $D$ be a simply connected domain with a Green's function. Write $D_t = F_t(D)$. Then
\begin{align*}
G_{D_t}(F_t(z_1),F_t(z_2)) - G_{D}(z_1,z_2)  = & \, \frac{4}{\pi} \Re \int_D \partial_zG_{D,z_1}(z)\partial_zG_{D,z_2}(z) \mu_t(z) d^2z 
+ o(|t|).
\end{align*}

If $f_t$ fixes $z_1$ and $\mu_t = 0$ in a neighborhood of $z_1$ for all sufficiently small $t$, then
\[
M_{D_t}(z_1)-M_{D}(z_1) = \frac{4}{\pi} \Re \int_D (\partial_z G_{D,z_1}(z) )^2 \mu_t(z) d^2 z +  o(|t|).
\]
\end{lemma}

In particular, if $d=\sum_{j=1}^na_j z_j$ is a divisor in $D$, and $f_t$ fixes neighborhoods of $z_1, \ldots, z_n$, then for $t$ sufficiently small,
\[
I(D_t, d) = I(D,d) + \frac{4}{\pi} \Re \int_D \left(\partial_z G_{D,d}(z) \right)^2 \mu_t(z) d^2 z + o(|t|).
\]

\section{Stationary solutions}

This section contains the proofs of Proposition~\ref{proB} and Theorem~\ref{ThmA} which together sets up the correspondence between stationary solutions and quadratic differentials. We prove the existence of stationary solutions by constructing a meromorphic quadratic differential which in turn is built by solving an extremal problem for the reduced Green's energy. We also study the particular half-translation surfaces that correspond to such quadratic differentials and prove Theorem~\ref{thm:1-1}.  
\subsection{Proof of Proposition~\ref{proB}}\label{sect:stationary-to-QD}
First note that the droplet components fill $\Sigma$ in the sense that their complement has measure $0$. Indeed, if this were not the case, there would exist some regular point of an interface with a prime end that does not correspond to a droplet component. This is impossible for a stationary solution since all source weights are strictly positive. Recall that we defined the quadratic differential $\vp$ as $\vp(z)dz^2:=\left(\partial_z G_{D^i,d^i}(z)\right)^2 dz^2$ on each droplet $D^i$. Assume now that $w$ is a regular interface point on $\partial D^i \cap \partial D^j$, and let $z$ be a local variable defined in a neighborhood $N$ of $w$. Since the solution is stationary we have that $$\partial_z G_{D^i,d^i}(z)dz=-\partial_z G_{D^j,d^j}(z)dz,$$ and hence $$\left(\partial_z G_{D^i,d^i}(z)\right)^2dz^2=\left(\partial_zG_{D^j,d^j}(z)\right)^2 dz^2$$
on the interface near $w$ which means that locally $\vp$ extends continuously over the interface. Since the interface is locally regular it follows e.g. from Morera's theorem that the local extension is holomorphic. 
Consider next an irregular interface point $w\in \Sigma$ where $k\geq 2$ droplets meet (the same droplet can be counted more than once). It is not hard to see, e.g., by considering the decay of harmonic measure near $w$, that stationarity and the assumption that droplet interfaces are piecewise regular implies that $k$ interfaces necessarily meet at the same angle at $w$, i.e., $2\pi/k$. From this it follows that irregular interface points correspond to points where individual droplet interfaces are irregular, and since these are assumed to be piecewise regular there are at most finitely many irregular interface points. By the angle bound $2\pi/k$ it follows that each $(\partial_z G_{D^i,d^i}(z))^2$ is bounded near an irregular point, and thus $\vp$ also extends holomorphically across these irregular points. Lastly, when $\Sigma$ is non-compact we also need to consider the punctures. If $k$ droplets meet at the puncture we still get that the interfaces meet at the same angle $2\pi /k$, and so if $k\geq 2$ we get as above that  $\vp$ extends holomorphically. If instead $k=1$ then by the piecewise regularity of the droplet interface one can bound $|(\partial_z G_{D^i,d^i})^2|$ locally showing that it extends meromorphically with a first order pole at the puncture.  
\qed
\subsection{Extremal problem and existence: Proof of Theorem~\ref{ThmA}}\label{sect:existenceproof}
Let $\Sigma$ be a Riemann surface of finite type. Fix an admissible droplet configuration $\mathcal{D}_0 = (D_0^i, \gamma_0^i, d^i)_{i=1}^n$ and recall the definition of reduced Green's energy in \eqref{def:green-energy-droplet-config}.  Let $\mathcal{H}(\mathcal{D}_0)$ denote the set of all admissible droplet configurations $(D^i,\gamma^i,d^i)$ such that for each $i$, $\gamma^i$ is homotopic to $\gamma^i_0$ in $\Sigma\smallsetminus \cup_i \textrm{supp}(d^i)$ (or in $\mathring{\Sigma}\smallsetminus \cup_i \textrm{supp}(d^i)$ if $\Sigma$ is non-compact).  We already noted that $I(\mathcal{D}_0) > - \infty$ but we also have the following upper bound.
\begin{lemma} \label{prop:bound1}
    For any collection of source divisors $(d^i)_{i=1}^n$ on $\Sigma$ (if $\Sigma = \CC$ or $\Sigma = \hat{\CC}$  we need to assume that $n>1$ and $n > 2$, respectively) one can find a constant $C<\infty$ depending only on $\Sigma$ and $(d^i)_{i=1}^n$ such that $I(\mathcal{D})\leq C$ for all admissible droplet configurations $\mathcal{D}=(D^i,\gamma^i,d^i)_{i=1}^n$. In fact each term $I(D^i, d^i)$ is uniformly bounded.
\end{lemma}

\begin{proof}
Let us first assume that $\Sigma$ is neither biholomorphic to $\mathbb{C}$ nor $\hat{\mathbb{C}}$.
It is enough to show that for any fixed pair of points $z_j,z_k\in \Sigma$ there is some constant $C$ such that for any simply connected domain $D\subset \Sigma$ we have that $M_D(z_j)\leq C$ and $G_D(z_j,z_k)\leq C$. We argue by contradiction. If $M_D(z_j)$ is not bounded then one could find a sequence of simply connected domains whose reduced moduli tended to $\infty$. By Montel's theorem this would give a holomorphic embedding of $\CC$ into $\Sigma$, which is impossible. If $G_D(z_j,z_k)$ is not bounded the same argument works.
If $\Sigma\cong \mathbb{C}$ we have the additional assumption of there being at least two source divisors. Let $w_2$ be a source point for $d^2$. Then we must have that $D^1\subseteq \Sigma\setminus\{w_2\}$. The argument above then yields a uniform bound on $I(D^1,d^1)$, and by symmetry we get a uniform bound on $I(\mathcal{D})$. 
If $\Sigma\cong \hat{\mathbb{C}}$ the additional assumption says that there are at least three source divisors. Let $w_2$ be a source point for $d^2$ and let $w_3$ be a source point for $d^3$. Then we must have that $D^1\subseteq \Sigma\setminus\{w_2,w_3\}$, and the argument works as before.
\end{proof}

We now turn to the proof of Proposition~\ref{thm2} which is the main step in the proof of Theorem~\ref{ThmA}. Before giving it we state and prove two topological lemmas that will be needed in the proof. 

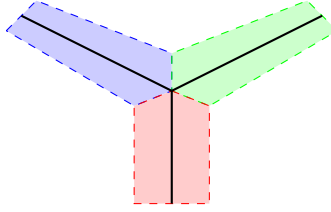
\begin{figure}

\begin{center}
    
\begin{tikzpicture}

\fill[blue!20!white] (-2.2,0.8) -- (-1.8,1.2) -- (0,0.5) -- (0,0) -- (-0.5,-0.2) -- (-2.2,0.8);

\fill[green!20!white] (2.2,0.8) -- (1.8,1.2) -- (0,0.5) -- (0,0) -- (0.5,-0.2) -- (2.2,0.8);

\fill[red!20!white] (0,0) -- (0.5,-0.2) -- (0.5,-1.5) -- (-0.5,-1.5) -- (-0.5,-0.2) -- (0,0);

\draw[blue, dashed] (-2.2,0.8) -- (-1.8,1.2) -- (0,0.5) -- (0,0) -- (-0.5,-0.2) -- (-2.2,0.8);

\draw[green, dashed] (2.2,0.8) -- (1.8,1.2) -- (0,0.5) -- (0,0) -- (0.5,-0.2) -- (2.2,0.8);

\draw[red, dashed] (0,0) -- (0.5,-0.2) -- (0.5,-1.5) -- (-0.5,-1.5) -- (-0.5,-0.2) -- (0,0);

\draw[thick] (0,0) -- (0,-1.5);

\draw[thick] (-2,1) -- (0,0);

\draw[thick] (0,0) -- (2,1);

\end{tikzpicture}

\end{center}
\caption{A tubular neighborhood of a graph near a vertex where three edges meet.} 
\label{fig:112}
\end{figure}

Let $\Gamma$ be a finite graph smoothly embedded in a differentialble surface $S$. 
\begin{df}
    A tubular neighborhood of an edge $\mathbf{e}$ is a homeomorphism between a neighborhood of $\mathbf{e}$ and a neighborhood of the unit square such that the image of $\mathbf{e}$ is the straight line between $(0,1/2)$ and $(1,1/2)$. With a tubular neighborhood $U$ of $\Gamma$ we mean a tubular neighborhood for each edge in $\Gamma$ such that whenever two or more edges meet at a vertex, the tubular neighborhoods glue linearly along the lines $\{0\}\times [0,1/2], \{0\}\times [1/2,1], \{1\}\times [0,1/2]$ and/or $\{1\}\times [1/2,1]$ (see Figure \ref{fig:112}). We identify $U$ with the union of the preimages of the unit squares. 
\end{df}

It is easy to see that such tubular neighborhoods always exist.

\begin{lemma} \label{lemcurve1}
    Let $\gamma_i:[0,1]\to [0,1]^2$ be finitely many continuous curves on the unit square with start and endpoints on the two vertical segments $\{0\}\times[0,1]$ and $\{1\}\times[0,1]$, and which do not intersect and do not self-intersect. Then one can find homotopic curves $\tilde{\gamma}_i$ with the same start and end points, that do not intersect and do not self-intersect, and such that for all $i$, if $\gamma_i$ ended at a different vertical segment than it started from, then $\tilde{\gamma}_i$ is a straight curve, while if $\gamma_i$ started and ended at the same vertical segment, then $\tilde{\gamma}_i$ is the concatenation of two straight curves. 
\end{lemma} 

\begin{proof}
    For simplicity let us first assume that $i=1,2$ and that $\gamma_i(0)=(0,a_i)$ and $\gamma_i(1)=(1,b_i)$. Without loss of generality we can assume that $a_1<a_2$ and that $a_1=b_1$. We can identify $[0,1]\times[0,1]$ with the vertical ends glued together with and annulus in $\mathbb{C}$, and then we see that $b_1<b_2$ follows from Jordan's curve theorem. If we now let $\tilde{\gamma}_i$ be the straight curve from $(0,a_i)$ to $(1,b_i)$ we see that the statements of the lemma hold in this case.
    
    If we instead have that $\gamma_i(0)=(0,a_i), \gamma_1(1)=(1,b_1)$ and $\gamma_2(1)=(0,b_2)$, then the same argument shows that $a_1<b_2$, and then we can let $\tilde{\gamma}_2$ go straight from $(0,b_1)$ to $(\epsilon,b_1)$ and then straight to $(0,b_2)$, where $\epsilon$ is chosen so that $\tilde{\gamma}_2$ do not intersect $\tilde{\gamma}_1$.

    Now assume that $\gamma_i(0)=(0,a_i)$ and $\gamma_i(1)=(0,b_i)$, and without loss of generality $a_1<a_2$. We then get closed simple curves $\gamma_i:[0,2]\to [-1,1]\times[0,1]$ by letting for $t\in[1,2]$, $\gamma_i(t):=(-\gamma_i(2-t)_1,\gamma_i(2-t)_2)$. Then by Jordan's curve theorem we cannot have $a_1<a_2<b_1<b_2$, and thus it is possible to construct curves $\tilde{\gamma}_i$ consisting of two straight segments that start and end at the given points and that do not intersect.

    The general case now follows from iterating the procedures above.
\end{proof}

\begin{lemma} \label{topcor}
    Let $\Gamma$ be a finite graph smoothly embedded in a differentiable surface $S$, let $U$ be a tubular neighborhoood of $\Gamma$, and let $\gamma_i$, $i=1,...,n$ be closed non self-intersecting curves in $U$ that do not intersect each other. Then there are closed curves $\tilde{\gamma}_i$ in $\Gamma\subset U$ that are non self-crossing and do not cross each other in $U$, that pass smoothly each edge they happen to enter, and such that each $\tilde{\gamma}_i$ is homotopic in $U$ to $\gamma_i$. 
\end{lemma}

\begin{proof}
First we note that from Lemma \ref{lemcurve1} we get non self-intersecting curves $\tilde{\gamma}_i$ in $U$ that do not intersect, such that for all i, $\tilde{\gamma}_i$ is homotopic to $\gamma_i$ in $U$, and restricted to a square of the tubular neighborhood of an edge the curve is either straight or composed of two straight segments. But this second case only happens when the curve comes back to the same vertical segment, and then that part can be retracted. If as a result of such a retraction we now get a curve which goes from one vertical side of  square to the other side and then back, then we retract that curve as well. Repeating this process we can assume that the restriction of each curve to a square of the tubular neighborhood of an edge the curve is straight.
Now we note that on each unit square of the tubular neighborhood we can homotope a curve $\gamma$ by letting $\gamma_s(t):=(\gamma(t)_1,(1-s)\gamma(t)_2+s/2)$. These homotopies glue together to give a homotopy $\gamma_s$ of any curve $\gamma$ in $U$ so that $\gamma_1$ lies on $\Gamma$. Applying this homotopy to the curves $\tilde{\gamma}_i$ now gives the Lemma, noting that the resulting curves are non-crossing as they are limits of non-intersecting curves.
\end{proof}

We are now ready to give the proof of Proposition~\ref{thm2}.
\begin{proof}[Proof of Proposition~\ref{thm2}]
We write $d=(d^i)_{i=1}^n$. Throughout, we fix coordinates near each point in the support of $d$. Set $M=\sup_{\mathcal{D} \in \mathcal{H}(\mathcal{D}_0)} I(\mathcal{D})$. By Lemma~\ref{prop:bound1}, $M < \infty$ and we also know that $M > -\infty$.
Let $\mathcal{D}_k \in \mathcal{H}(\mathcal{D}_0)$ be a sequence of droplet configurations 
such that $\lim_{k \to \infty} I(\mathcal{D}_k) = M$. Let $1 \le i \le n$ and fix some point $p \in \textrm{supp}(d^i)$. For each $k$, let $E_{w,k}=E_{w,k}^i$ be the droplet component of $\mathcal{D}_k$ which contains $p$. By Lemma~\ref{prop:bound1} each of the terms in $I(\mathcal{D}_k)$ is uniformly bounded above. Hence the conformal radius of $E_{p,k}$ seen from $p$ is uniformly bounded away from $0$ and $\infty$ as $k \to \infty$. Next, let $f_{p,k} : \mathbb{D} \to E_{p,k}$ be the conformal map taking $0$ to $p$ with positive derivative there. Using Lemma~\ref{prop:bound1} with Montel's theorem and the Carath\'eodory kernel theorem, we can find a subsequence of the maps $f_{p,k}$ which converges locally uniformly on $\mathbb{D}$ to a conformal map $f_p=f_{p, \infty}$, as $k \to \infty$. This limiting $f_p$ maps $\mathbb{D}$ onto a simply connected hyperbolic domain $E_p$ which contains $p$. We obtain a divisor $e_p$ by considering in addition to $p$ those points in $d^i$ lying in $E_p$.  

Repeating this process for the remaining points in $\textrm{supp}(d^i)$ and for all $i$, taking further subsequences if necessary, results in a finite set of simply connected domains $(E^{j})_{j=1}^m$ with divisors $(e^j)_{j=1}^m$ (keeping the weights as in $d$) each supported in $E^j$. Note that $\cup_j\textrm{supp}(e^j) = \cup_i\textrm{supp}(d^i)$. Using the continuity properties of the reduced modulus and Green's function with respect to Carath\'eodory convergence, we see that  
\[I(E,e) = \sum_{j=1}^mI(E^j, e^j)= M.\] 
We will now carry out a quasiconformal variation of $\Sigma$. Let $N \subset \Sigma$ be a simply connected set, chosen so small that it is contained in some coordinate patch, which does not intersect $\textrm{supp}(d)$. Choose a coordinate $z:N \to z(N)$ and consider for $h \in C^1_0(N)$ and small $t$, the map $f_t:\Sigma \to \Sigma$
\[
f_t(z) = z + t h(z),
\]
extended to the identity outside of $N$.
For small enough $t$, $f_t$ is a quasiconformal homeomorphism of $\Sigma$ and the corresponding Beltrami coefficient is
\[
\mu_{t}(z)= \frac{t \partial_{\bar z} h(z)}{1+ t \partial_{z} h(z)}.
\]
For small $t$, $f_t$ deforms the domains $E^j, j=1,\ldots,m,$ homeomorphically with some neighborhood of each point in $\textrm{supp}(e)$ fixed. Suppose for some fixed $j$, $N \cap E^j \neq \emptyset$ and let $w: E^j \to \mathbb{D}$ be a conformal map sending $z^j_1$ to $0$. By considering the inverse of this map composed with the homeomorphism $N \to z(N)$, we obtain a conformal map $z(w): w(N) \to z(N)$. 
Write $E_t^j = f_t(E^j)$ and $w_t: E_t^j \to \mathbb{D}$ for the corresponding conformal map sending $z^j_1$ to $0$. Then $F_t = w_t \circ f_t \circ w^{-1}$ is a quasiconformal map  $\mathbb{D} \mapsto \mathbb{D}$ with Beltrami coefficient given by\[
 \mu_{F_t}(w) = \mu_{t}(z(w))\frac{ \overline{z'(w)} }{z'(w)}.
\]
(The dilatation is not affected by postcomposing by the conformal map $w_t$.)
Let $z_1,z_2 \in \textrm{supp}(e^j)$. Then $f_t(z_1)=z_1, f_t(z_2)=z_2$, so
\[
G_{E_t^j}(z_1,z_2) - G_{E^j}(z_1,z_2) = G_{\mathbb{D}}(F_t(w_1),F_t(w_2)) - G_{\mathbb{D}}(w_1,w_2),
\]
where $w_1=w(z_1), w_2=w(z_2)$.
On the other hand, by Lemma~\ref{lem:variation}, as $t \to 0$,
\begin{align*}
G_{\mathbb{D}}(F_t(w_1),F_t(w_2)) - G_{\mathbb{D}}(w_1,w_2) &= \frac{4}{\pi} \Re \int_{\mathbb{D}}\partial_wG_{\mathbb{D}}(w,w_1)\partial_wG_{\mathbb{D}}(w,w_2) \mu_{F_t}(w) dA(w) + O(|t|^2). 
\end{align*}
Similarly, using that $\mu_{F_t} = 0$ in a neighborhood of $w_1$ so that $F_t$ is conformal at $w_1$ (and the covariance terms cancel),
\[
M_{E^j_t}(z_1) - M_{E^j}(z_1) =\frac{4}{\pi}\Re \int_{\mathbb{D}}(\partial_wG_{\mathbb{D}}(w,w_1)^2 \mu_{F_t}(w) dA(w) + o(|t|). 
\]
It follows that
\begin{align}
I(E_t^j, e^j) - I(E^j,e^j)  & = \frac{4}{\pi} \Re \int_{\mathbb{D}}\left(\sum_k a^i_k \partial_w G_{\mathbb{D}}(w,w_k^j) \right)^2\mu_{F_t}(w) dA(w) + O(|t|^2) \nonumber \\
& = \frac{4}{\pi} \Re \int_{\mathbb{D}}\left(\sum_k a^j_k\partial_w G_{\mathbb{D}}(w,w_k^j) \right)^2\mu_{t}(z(w))\frac{\overline{z'(w)}}{z'(w)} dA(w)  + O(|t|^2) \nonumber \\
&=  \frac{4}{\pi} \Re \int_{\mathbb{D}}\left(\sum_k a^j_k \partial_z G_{E^j}(z(w),z_k^j) \right)^2\mu_{t}(z(w))|z'(w)|^2 dA(w) + O(|t|^2) \nonumber\\
& =    \frac{4t}{\pi} \Re \int_{z(N\cap E^i)}\left(\partial_z G_{E^j,e^j}(z) \right)^2\partial_{\bar z}h(z) dA(z)+ O(|t|^2). \label{eq:variation1}
\end{align}
Define a quadratic differential on $\Sigma$ by setting
\[
\varphi(z)dz^2 :=  \left(\partial_z G_{E^i,e^i}(z) \right)^2dz^2, \quad z \in E^j, \quad j=1,\ldots,m,
\]
and $\varphi(z) := 0, z \in \Sigma \smallsetminus (\cup_{j=1}^m E^j)$. Then $\varphi(z) dz^2$ is meromorphic in each $E^j$ with second order poles at $\textrm{supp}(e^j)$.
Summing over $j$ and using \eqref{eq:variation1}, we obtain
\begin{equation}\label{eq:weyl1}
 I(E_t, e) - I(E,e) =   \frac{4 t}{\pi} \Re \int_{N}\varphi \partial_{\bar z}h dA + O(|t|^2).
\end{equation}
By the extremal property of $(E,e)$ we have $I(E_t, e) - I(E,e) \le 0$ and since $h$ was an arbitrary $C^1_0(N)$ function, \eqref{eq:weyl1} implies that
\begin{equation}\label{eq:weyl2}
 \int_{N}\varphi \partial_{\bar z}h dA = 0.
\end{equation}
 The set $N$ was also arbitrary except for not intersecting $\textrm{supp}(e)$, so \eqref{eq:weyl2} and Weyl's lemma (see, e.g., \cite[Chapter 10.3]{gardiner}) imply that there exists a quadratic differential $\tilde \varphi$ which is meromorphic on $\Sigma$ and a.e. equal to $\varphi$. It follows that $\tilde \varphi dz^2 = (\partial_z G_{E^j,e^j})^2 dz^2$ in $E^j$ and $\Sigma \setminus \cup_{j=1}^m E^j$ has measure $0$. 
 
 It now only remains to show that the $E^j$ are the simply connected components of an admissible droplet configuration $\mathcal{D}_\infty \in \mathcal{H}(\mathcal{D}_0)$. 
 For this we will use Lemma~\ref{topcor}. 

First we assume that $\Sigma$ is closed. Let $\Gamma:=\Sigma\setminus \cup_j E^j$. Then $\Gamma$ is the critical graph of a meromorphic quadratic differential on $\Sigma$ which in particular is a smoothly embedded finite graph. Pick a tubular neighborhood $U$ of $\Gamma$ which does not contain any point in $\textrm{supp}(d)$. From the proof of the existence of the quadratic differential we see that for $k$ large enough (along a subsequence) each droplet interface $\gamma^i_k$ will lie in $U$. Fix such a $k$ and consider the curves $\gamma^i_k$. By Lemma~\ref{topcor} we get curves $\tilde{\gamma}^i$ in $\Gamma$ that are homotopic to $\gamma^i_k$ in $U$ and hence also in $\Sigma\setminus \textrm{supp}(d)$, and from Lemma~\ref{topcor} it also follows that $\mathcal{D}_\infty:=(D^i,\tilde{\gamma}^i,d^i)$ will be an admissible droplet configuration. 

Now we consider the case when $\Sigma$ is non-compact. Let $\Gamma:=\tilde{\Sigma}\setminus \cup_j E^j$. Then $\Gamma$ is the critical graph of a meromorphic quadratic differential on $\tilde{\Sigma}$ which in particular is a smoothly embedded finite graph. We first define an embedded graph $\overline{\Gamma}$ on $\mathring{\Sigma}$ as being equal to $\Gamma$ on $\Sigma$ and if $\Gamma$ contains a puncture $p_i$ then we add to $\overline{\Gamma}$ the corresponding boundary circle in $\mathring{\Sigma}$. Now we note that $\overline{\Gamma}$ is a smoothly embedded graph in the differentiable surface (without boundary) $\Sigma_+\supset \mathring{\Sigma}$ that we defined in Section \ref{subsec:prel}. We pick a tubular neighborhood $U$ of $\overline{\Gamma}$ in $\Sigma_+$ which does not contain any point in $\textrm{supp}(d)$. As before, for $k$ large enough (along a subsequence) each droplet interface $\gamma^i_k$ will lie in $U$. Fix such a $k$ and consider the curves $\gamma^i_k$. By Lemma~\ref{topcor} we get curves $\tilde{\gamma}^i$ in $\overline{\Gamma}$ that are homotopic to $\gamma^i_k$ in $U$ and hence also in $\mathring{\Sigma}\setminus \textrm{supp}(d)$, and from Lemma~\ref{topcor} it again follows that $\mathcal{D}_\infty:=(D^i,\tilde{\gamma}^i,d^i)$ will be an admissible droplet configuration.
\end{proof}

\begin{proof}[Proof of Theorem~\ref{ThmA}]
Let $\mathcal{D}_\infty \in \mathcal{H}(\mathcal{D}_0)$ be as in Proposition~\ref{thm2}. In each droplet component $E^j$ we have \[\varphi(z) dz^2 = (\partial_z G_{E^j,d^j}(z))^2 dz^2\] and this determines a meromorphic quadratic differential on $\Sigma$. If $z$ is a regular point on a droplet interface then it is not a critical point of $\varphi$ and hence $\varphi$ is holomorphic at $z$. This immediately implies that $\mathcal{D}_\infty$ is a stationary solution.
 \end{proof}

\subsection{Half-translation surfaces of Green's type} \label{subsec:half}

A quadratic differential on a Riemann surface yields a representation of that Riemann surface as a half-translation surface, see, e.g., \cite{Zor}. In this section we will describe directly those half-translation surfaces that correspond to a quadratic differential $\vp$ associated to a stationary solution.

Let $(D^i,\gamma^i)$ be a stationary solution on $\Sigma$ with source divisors $d=(d^i)$. 
\begin{prop}
Assume that each source divisor $d^i$ is a singleton divisor $d^i=a^iz^i$. Then the corresponding half-translation surface will consist of half-strips $S^i:=\{ \Re(z)<0, 0\leq \Im(z)\leq \pi a^i\}$ with the top and bottom of each half strip being identified by the translation $z\mapsto z+\pi ia^i$, and with the vertical boundary of $S^i$ corresponding to the droplet interface $\gamma^i$. When two interfaces $\gamma^i$ and $\gamma^j$ intersect along an arc this corresponds to a piece of the vertical boundary of $S^i$ being identified via a minus-translation to a piece of the vertical boundary of $S^j$. 
\end{prop}

\begin{proof}Let $w$ be an conformal map from $D^i$ to the unit disc $\mathbb{D}$ mapping $z_i$ to the origin. It is then clear that in this coordinate $G_{D^i,d^i}(w)=-a^i\log|w|$. If $z=x+iy$ is the coordinate on $S^i$ then $e^{2z/a^i}$ maps $S^i$ to $\mathbb{D}^*$, and in this coordinate we get that $$G_{D^i,d^i}(z)=-a^i\log|e^{2z/a^i}|=-2x$$ and hence $$(\partial_z G_{D^i,d^i}(z))^2dz^2=dz^2.$$

It is clear that the vertical boundary of $S^i$ corresponds to the droplet interface $\gamma^i$. The fact that we have a half-translation surface means that whenever two half-strips $S^i$ and $S^j$ are identified along their vertical boundaries this must be via a minus-translation. (see Figure \ref{fig:5} in Section \ref{subsec:st}).
\end{proof}

We now consider the case when not all source divisors are singletons. Assume that the droplet $D^i$ has several source points $z^i_j$, $j=1, \ldots, n_i$, with weights $a^i_j$. We claim that $D^i$ can be represented as a partial translation surface\footnote{Recall from Section \ref{subsec:prel} that a partial translation surface is given by a collection of polygons with some but not all sides being pairwise identified via translations.} $T$ of a simple form so that $dz$ on $T$ corresponds to $-\partial_z G_{D^i,d^i}(z)dz$ on $D^i$.

The class of partial translation surfaces we will consider can be described as follows. 
\begin{df}[Green's surface] The class of Green's surfaces is the smallest class of partial translation surfaces satisfying the following conditions. Half-strips $\{z: \Re(z)\leq 0, 0\leq  \Im(z)\leq \pi a\}$ with the identification $z=z+\pi i a$ belong to the class. In general, any element in the class has a single unpaired side and we call this the \emph{boundary side}. The boundary side is supposed to be vertical and to have the surface to the left. For any finite collection of surfaces $T^i$ in the class with boundary side lengths $\pi a_i$ and a rectangle $$R:=\{z: -b\leq \Re(z)\leq 0, 0\leq \Im(z)\leq \pi a\}$$ such that $a=\sum_i a^i$ and with the top and bottom sides of $R$ being identified by $z=z+\pi i a$, then the partial translation surface obtained by gluing the surfaces $T^i$ to the rectangle along its left side $\Re(z)=-b$ (in any chosen order) also belongs to the class. 
\end{df}

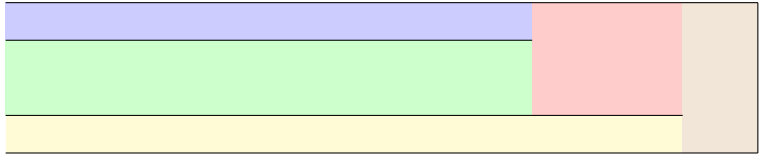
\begin{figure}

\begin{center}
    
\begin{tikzpicture}

\fill[blue!20!white] (-7,1) -- (0,1) -- (0,0.5) -- (-7,0.5) -- (-7,1);

\fill[green!20!white] (-7,0.5) -- (0,0.5) -- (0,-0.5) -- (-7,-0.5) -- (-7,0.5);

\fill[red!20!white] (0,1) -- (2,1) -- (2,-0.5) -- (0,-0.5) -- (0,1);

\fill[yellow!20!white] (-7,-0.5) -- (2,-0.5) -- (2,-1) -- (-7,-1) -- (-7,-0.5);

\fill[brown!20!white] (2,1) -- (3,1) -- (3,-1) -- (2,-1) -- (2,1);

\draw (-7,1) -- (3,1);

\draw (-7,0.5) -- (0,0.5);

\draw (-7,-0.5) -- (2,-0.5);

\draw (-7,-1) -- (3,-1);

\draw (3,1) -- (3,-1);

\end{tikzpicture}

\end{center}
\caption{A representation of a Green's surface.} 
\label{fig:111}
\end{figure}

Figure \ref{fig:111} shows a representation of a Green's surface where the blue, green and yellow parts represent half-strips while the red and brown parts are rectangles.

\subsubsection*{Green's functions on Green's surfaces} On any polygonal piece of a (partial) translation surfaces one can specify an $x$-coordinate, but this function is only determined up to the addition of a real constant, and it is not always possible to patch these coordinates to get a well-defined $x$-coordinate on the whole surface. On a Green's surface $T$ however there is a canonical choice of $x$-coordinate, that we will write $x_T$, and we will call the $-2x_T$ the \emph{Green's function} on the Green's surface. If the surface $T$ is a half-strip, then we let $x_T$ be the $x$-coordinate which is zero on the boundary side. If the surface $T$ consists of finitely many Green's surfaces $T^i$ glued to the left side of the rectangle $R:=\{z: -b\leq \Re(z)\leq 0, 0\leq \Im(z)\leq \pi a\}$, and we assume that we have already defined $x_{T_i}$ for each $T^i$, then we define $x_T$ to be the $x$-coordinate on $R$ which is zero on its right side, while on each $T_i$ we define $x_T$ as $x_T:=x_{T_i}-b$. In this way we can define the canonical $x$-coordinate and the corresponding Green's function on all Green's surfaces.

\subsubsection*{Representing a droplet as a Green's surface} To see that a droplet $D^i$ can be represented as a Green's surface so that $-\partial_z G_{D^i,d^i}(z)dz$ becomes $dz$, we argue by induction over the number of source points. The case of only one source point was shown above. Assume that it is known for up to $k$ source points, and that $D^i$ has $k+1$ source points. Let $A_x:=\{w\in D^i: G_{D^i,d^i}>-2x\}$ and let $x_0$ be infimum of $x$ so that $A_x$ only has one connected component. Thus $A_{x_0}$ will have several components $A^1,\ldots, A^n$. Each component will contain at least one source point and it follows from the maximum principle that all components $A^j$ will be simply connected. Also, on $A^j$ we will have that $$G_{D^i,d^i}+2x=G_{A^j,d^i_j},$$ where $d^i_j$ is the restriction of the source divisor to $A^j$. By the induction hypothesis each $A^j$ can be represented by a translation surface $T^j$ in the class. On the other hand we can represent $D^i\smallsetminus \overline{A_{x_0}}$ by the rectangle $R=\{z: x\leq\Re(z)\leq0, 0\leq\Im(z)\leq\pi a\}$ with he identification $z=z+\pi i a$ and $a$ being the sum of the source strengths. One thus sees that $D^i$ can be represented by the surface one gets by gluing the surfaces $T^j$ to the rectangle $R$ along the level set $G_{D^i,d^i}=-2x$.

\subsubsection*{A one-to-one correspondence} Finally we can represent $\Sigma$ as a half-translation surface by taking the Green's surfaces for each $D^i$ and then glue them in the same way as when there only was one source point in each droplet.
\begin{df}[Half-translation surfaces of Green's type]We say that a half-translation surface consisting of Green's surfaces with minus-translation identifications along their vertical right-most boundaries is a half-translation surfaces of \emph{Green's type}.  
\end{df}

We are ready to prove Theorem~\ref{thm:1-1}.

\begin{proof}[Proof of Theorem~\ref{thm:1-1}] We need to show that if one can represent a Riemann surface $\Sigma$ as a half-translation surface $T$ of Green's type, then this corresponds to a stationary solution to the competitive Hele-Shaw flow on $\Sigma$. To see this assume that $T$ consists of Green's surfaces $T^i$ glued together along their right-most vertical boundaries, and let $D^i$ denote the corresponding domains in $\Sigma$. Each $T^i$ will contain a finite number of half-strips with height $\pi a^i_l$, and the half-strips will in $\Sigma$ correspond to discs punctured at points $z^i_j$. One can now easily check that the Green's function of $T^i$ is equal to $G_{D^i,d^i}$ where $d^i:=\sum_j a^i_jz^i_j$, and from this it follows immediately that $(D^i,\gamma^i)$ is a stationary solution to the competitive Hele-Shaw flow on $\Sigma$ with driving divisor $d=(d^i)$.
\end{proof}

\section{Special cases and examples}
\subsection{Two droplets in $\hat{\mathbb{C}}$ separated by a Jordan curve}\label{sect:jordan}
In this section we will study the special case of stationary Jordan curves, i.e., solutions with two droplets in $\hat{\mathbb{C}}$, separated by a Jordan curve $\gamma$ in $\mathbb{C}$. In this case other descriptions are easy to obtain.

 Write $D, D^*$, respectively, for the bounded and unbounded component of $\mathbb{C} \smallsetminus \gamma$, and let two corresponding divisors $d = \sum_j a_j z_j, d^* = \sum_j b_j w_j$ be given. We assume that $\infty$ is in the support of $d^*$ and $0$ is in the support of $d$ and write $a_1,b_1$ for the weights at $0$ and $\infty$, respectively.
We define
\[
I(\gamma, d,d^*) = I(D, d) + I(D^*, d^*).
\]

The simplest example is when $d=1 \cdot 0, d^* = 1 \cdot \infty$. Let $f:\mathbb{D} \to D, \, g:\mathbb{D}^* \to D^*$ be conformal maps fixing $0$ and $\infty$. Then  
\[
I(\gamma, d, d^*) = M_D(0) + M_{D^*}(\infty) = \log r_D(0) + \log r_{D^*}(\infty) \le 0
\]
with equality if and only if $\gamma$ is any circle centered at $0$, a consequence of the Grunsky inequality. This case is somewhat degenerate in the sense that any circle gives energy $0$. The corresponding quadratic differential is $\varphi(z) = z^{-2}dz^2$.
\subsubsection*{Variational formulas and energy monotonicity}
Before studying stationary solutions we will consider Hadamard variations corresponding to the competitive Hele-Shaw dynamics. Such a variation is well-defined assuming the curve is sufficiently smooth without information about local existence of solutions to the competitive Hele-Shaw problem itself. Indeed, only sufficient smoothness of the relevant vector field along the curve is needed for the Hadamard variation to be well-defined.

We follow \cite{GV} and use the notation $a \nabla_{D,z}$ for the Hadamard variation corresponding to Hele-Shaw injection in $D$ with weight $a$ at $z$. That is, the normal velocity at $\zeta \in \gamma=\partial D$ in the outward pointing direction is given by $a$ times the Poisson kernel, so that in the notation of Lemma~\ref{lemma:hadamard}, $\nu(\zeta) = aP_{D,z}(\zeta)$. Any source divisor $d= \sum_{j=1}^n a_j z_j$ then gives rise to a Hadamard variation by setting
\[
\nabla_{D,d}:=\sum_{j=1}^n a_j \nabla_{D,z_j}.
\]
Recall that the reduced Green's energies $I(D,d)$ and $I(D^*,d^*)$ are increasing in $D$ and $D^*$, respectively. The next result shows that the competitive Hele-Shaw variation of their sum $I(\gamma, d,d^*)$ is nevertheless non-negative. This observation motivated the proof of Theorem~\ref{ThmA}.

\begin{prop}\label{prop:gradient}
Let $\gamma$ be a smooth Jordan curve and let $d,d^*$ be divisors in $D, D^*$, respectively. Then
\[
\left( \nabla_{D,d} + \nabla_{D^*,d^*} \right) I(\gamma, d, d^*) \ge 0,
\]
with equality if and only if $\gamma$ is a stationary solution for the competitive Hele-Shaw problem.
\end{prop}
\begin{proof}
If $\gamma$ is smooth, then the Poisson kernel is also smooth so the Hele-Shaw vector fields from both sides are smooth. Hadamard's formula for the variation of the Green's function can therefore be applied.
For $j=1,\ldots, n$, set $P_j(\zeta):=P_{D,z_j}(\zeta)$ and for $j=1,\ldots,m$ set $Q_j(\zeta) = P_{D^*,w_j}(\zeta)$ for the Poisson kernels. We get
\begin{align*}
 \left( \nabla_{D,d} + \nabla_{D^*,d^*} \right) I(\gamma, d, d^*)  &= \left(\sum_{j=1}^n a_j \nabla_{D,z_j} + \sum_{j=1}^m b_j \nabla_{D^*,w_j} \right) \left( I(D, d) + I(D^*, d^*) \right) \\
 &= \frac{1}{2\pi} \int_\gamma\left[ \left(\sum_{j=1}^n a_j P_j\right)^2- \left(\sum_{j=1}^m b_j Q_j\right)^2 \right]\left(\sum_{j=1}^na_j P_j - \sum_{j=1}^mb_j Q_j   \right) |d\zeta| \\
  &= \frac{1}{2\pi} \int_\gamma \left( \sum_{j=1}^n a_j P_j + \sum_{j=1}^m b_j Q_j \right)\left(\sum_{j=1}^na_j P_j - \sum_{j=1}^mb_j Q_j   \right)^2 |d\zeta|   \ge 0.
\end{align*}
Equality holds if and only if $\zeta \mapsto \sum_{j=1}^na_j P_j(\zeta)- \sum_{j=1}^mb_j Q_j(\zeta) \equiv 0$ along $\gamma$, which means that precisely that $\gamma$ is a stationary solution.
\end{proof}

We can easily compute the competitive Hele-Shaw variations of the basic geometric functionals. 
\begin{prop}\label{prop:area-variation}
Let $\gamma$ be a smooth Jordan curve and let $d,d^*$ be divisors in $D, D^*$, respectively. Write $A$ for the area of $D$ and $L$ for the perimeter, $|\partial D|$. Then
\[
\left( \nabla_{D,d} + \nabla_{D^*,d^*} \right) A = 2\pi(|d|-|d^*|)
\]
and
\[
\left( \nabla_{D,d} + \nabla_{D^*,d^*} \right) L = 2\pi(\sum_{j=1}^n a_j\kappa_D(z_j)- \sum_{j=1}^m b_j\kappa_{D^*}(w_j)),
\]
where $\kappa_D(z), \kappa_{D^*}(z)$ are the harmonic extensions of the curvature of $\partial D$ into $D$ and $D^*$, respectively. 
\end{prop}
\begin{proof}
We have $A = \int_{D} dx dy$ and since the Poisson kernels integrate to $2\pi$, 
\[
\left( \nabla_{D,d} + \nabla_{D^*,d^*} \right)  A = \int_{\gamma} \left(\sum_{j=1}^n a_jP_{D,z_j}(\zeta) - \sum_{j=1}^m b_jP_{D^*, w^j}(\zeta) \right) |d\zeta| = 2\pi(|d|-|d^*|).
\]

Next, we have $L = \int_{\partial D} |d\zeta|$, so if $\kappa(\zeta)$ denotes the geodesic curvature of $\gamma$ at $\zeta$, then the well-known formula for the variation of the arc-length element gives
\[
\left( \nabla_{D,d} + \nabla_{D^*,d^*} \right) L =  \int_{\gamma} \kappa(\zeta)\left(\sum_{j=1}^n a_jP_{D,z_j}(\zeta) - \sum_{j=1}^m b_jP_{D^*, w^j}(\zeta) \right) |d\zeta|
\]
which gives the stated formula.
\end{proof}
Proposition~\ref{prop:area-variation} implies any solution to competitive Hele-Shaw problem is locally area preserving if $|d| = |d^*|$. In addition, we see that $|d|=|d^*|$ is a necessary condition for the existence of a stationary solution in this setting.

\subsubsection*{Level sets and lemniscates}
We now give several alternative descriptions of stationary Jordan curves. While our point of view is new, most of the results here can be seen to follow easily from work of Younsi \cite{younsi}. Recall that we assume throughout that $\infty$ with weight $b_1$ is in the support of $d^*$.

\begin{prop}\label{prop:level-line}
    Suppose the Jordan curve $\gamma$ is a stationary solution to the competitive Hele-Shaw problem with divisors $d,d^*$. Then $\gamma$ is a level line of the potential
    \[
    R(z) = \sum_{j=2}^m b_j \log|z-w_j|^{-1}-\sum_{j=1}^n a_j \log|z-z_j|^{-1}.
    \]
  Conversely, if the Jordan curve $\gamma$ is a level line of $R$ separating $\textrm{supp}(d)$ from $\textrm{supp}(d^*)$ then $\gamma$ is a stationary solution to the competitive Hele-Shaw problem with divisors $d,d^*$.   
\end{prop}
\begin{rem}
    In the special case when $d=\sum_{j=1}^n n^{-1} \cdot z_j$ and $d^* = 1 \cdot \infty$, it follows from the proposition that any stationary $\gamma$ is a polynomial lemniscate. Similarly, if all weights are rational, $\gamma$ is a rational lemniscate. See \cite{younsi} and the references therein.
\end{rem}
\begin{rem}
By the second part of the proposition, a ``typical'' stationary solution $\{z : R(z) = \lambda_0 \}$ gives rise to a family of solutions $\lambda \mapsto \gamma_\lambda:=\{z: R(z) = \lambda \}$, by varying the level $\lambda$ in a neighborhood of $\lambda_0$. In particular, stationary solutions are not unique in this case. We expect uniqueness however if the area of the droplet is prescribed.
\end{rem}
\begin{proof}[Proof of Proposition~\ref{prop:level-line}]By assumption the Jordan curve $\gamma$ is a stationary solution and hence an analytic curve. Therefore, by Proposition~\ref{prop:area-variation}, we have $d=d^*$. Define the function $G(z) :=  G_{D^*, d^*}(z), z \in D^*$ and $G(z) := - G_{D, d}(z), z \in \overline{D}$. Set
\[    R(z) = \sum_{j=2}^m b_j \log|z-w_j|^{-1}-\sum_{j=1}^n a_j \log|z-z_j|^{-1}.
\]
Since $d=d^*$, \[R(z) = \left(\sum_{j=1}^na_j-\sum_{j=2}^mb_j\right) \log|z| + O(1) = b_1 \log|z| + O(1),\] as $|z| \to \infty$. Moreover, since $\gamma$ is stationary, $\partial_z (G-R)$ extends to a holomorphic function in $\hat{\mathbb{C}}$. So  
\[H(z) := G(z) -R(z)\] extends to a bounded harmonic function in in $\hat{\mathbb{C}}$, and hence $H$ is constant. On the other hand, $G(z)=0$ along $\gamma$, so it follows that $R(z)$ is constant along $\gamma$, as claimed.
 Suppose now that the Jordan curve $\gamma=\{z : R(z) = \lambda \}$ is a level line of $R$ separating $\textrm{supp}(d)$ from $\textrm{supp}(d^*)$ and write $D, D^*$ for the bounded and unbounded component of $\hat{\mathbb{C}} \smallsetminus \gamma$, respectively. Then $\gamma$ is smooth and by uniqueness of the Green's function it follows that \[G_{D^*,d^*}(z) = (R(z)-\lambda)\mid_{D^*}, \quad G_{D,d}(z) = (\lambda-R(z))\mid_{D}.\]
 Hence $\partial_zG_{D,d} = -\partial_zG_{D^*,d^*}$ along $\gamma$, which means it is a stationary solution.
\end{proof}

Given a Jordan curve $\gamma$, the welding homeomorphism $h : S^1 \to S^1$ associated to it is defined by $h=g^{-1} \circ f\mid_{S^1}$, where $f:\mathbb{D} \to D$ and $g:\mathbb{D}^* \to D^*$ are conformal maps as above. Conversely, given an orientation preserving homeomorphism $h:S^1 \to S^1$, the conformal welding problem is to find a Jordan curve and corresponding conformal maps so that $h=g^{-1} \circ f\mid_{S^1}$ holds. 

Let us make some remarks on the welding homeomorphism corresponding to a given stationary curve $\gamma$ as in the beginning of the section. 

Write $x_j=f^{-1}(z_j), \, y_j = g^{-1}(w_j)$. Then stationarity and conformal covariance of the Poisson kernel immediately gives the following functional equation for $h$.
\begin{lemma}\label{lem:welding-functional-eq}
We have
\[
\sum_{j=1}^ma_j  P_{\mathbb{D}, x_j}(z) = \sum_{j=1}^nb_j  |h'(z)|P_{\mathbb{D}^*, y_j}(h(z)), \quad z \in \partial \mathbb{D}.
\]
\end{lemma}
Various relations can be derived from this. We give one example. For more, see \cite{younsi}.
\begin{prop}
Suppose that $d = \sum_{j=1}^m a_j z_j$ and $d^* = 1 \cdot \infty$. Suppose $\gamma$ is a stationary Jordan curve with divisors $d,d^*$ and let $h$ be the welding homeomorphism of $\gamma$. Then
    \[
h(z)= \lambda \prod_{j=1}^m\left(\frac{z - x_j}{1-\overline{x_j}z} \right)^{a_j}, 
\]
where for $j=1,2,\ldots m$, $x_j = f^{-1}(z_j)$ and $|\lambda|=1$.
\end{prop}
\begin{proof}
Set $M_j(z) = (z-x_j)/(1-\overline{x_j}z)$. Then
\[
P_{\mathbb{D}, x_j}(e^{i\theta})=\frac{d}{d\theta} \log M_j(e^{i\theta}). 
\]
Write $h(e^{i\theta})=:e^{i\psi(\theta) + i\psi_0}$ with $\psi: [0,2\pi) \to [0,2\pi)$  an increasing homeomorphism. Since $ie^{i\theta} h'(e^{i\theta}) = \frac{d}{d\theta} h(e^{i\theta}) = i \psi'(\theta) h(e^{i\theta})$, it follows that $|h'(e^{i\theta})| = |\psi'(\theta)| = \psi'(\theta)$, and hence Lemma~\ref{lem:welding-functional-eq} gives
\[
\sum_{j=1}^ma_j  \frac{d}{d\theta} \log M_j(e^{i\theta}) =  \psi'(\theta).
\]
Integrating this gives the claim.
\end{proof}

\subsection{Examples with three droplets in $\hat{\mathbb{C}}$}

Above we have seen explicit examples with two droplets in $\hat{\mathbb{C}}$. We now look at examples with three droplets in $\hat{\CC}$. 

We assume that all three source divisors are singletons. Then, without loss of generality we can assume that the three source divisors are $d^0=a_0\cdot 0$, $d^1=a_1\cdot 1$ and $d^2=a_{\infty}\cdot\infty$. It is easy to check that the only quadratic differential $\varphi$ on $\hat{\CC}$ which is holomorphic except having double poles at $0,1,\infty$ and such that $\textrm{Res}_{0}(\sqrt{\vp})=\pm a_0/2, \textrm{Res}_{1}(\sqrt{\vp})=\pm a_1/2$ and $\textrm{Res}_{\infty}(\sqrt{\vp})=\pm a_{\infty}/2$ is given by 
\begin{equation} \label{eq:fi}
\vp(z)=\frac{a_0^2+(a_1^2-a_0^2-a_{\infty}^2)z+a_{\infty}^2z^2}{4z^2(z-1)^2}dz^2.
\end{equation}

Thus if $a_0^2+(a_1^2-a_0^2-a_{\infty}^2)z+a_{\infty}^2z^2$ has zeroes $z_i$ then three or more droplets will meet at exactly these points. To trace the interfaces between these meeting points we follow the vertical foliation by solving the equation 
\begin{eqnarray*}
\textrm{Re}\left(\int_{z_i}^w\frac{\sqrt{a_0^2+(a_1^2-a_0^2-a_{\infty}^2)z+a_{\infty}^2z^2}dz}{2z(z-1)}\right)=0.
\end{eqnarray*} 

\subsection{Examples with two droplets in $\mathbb{C}$} \label{sec:ex}

If we let $a_{\infty}$ tend to zero in the formula (\ref{eq:fi})  we see that $$\vp(z)=\frac{a_0^2+(a_1^2-a_0^2)z}{4z^2(z-1)^2}dz^2$$
describes a stationary solution $(D^i,\gamma^i)$ of the competitive Hele-Shaw problem in $\CC$ with the two source divisors $d^0=a_0\cdot 0$ and $d^1=a_1\cdot 1$. 

Let us look more closely at the case $a_0=2$ and $a_1=1$. We see that three droplets will meet at the zero $z_0=4/3$. To find the interfaces emanating from this point we solve the equation $$\textrm{Re}\left(\int_{4/3}^w\frac{\sqrt{4-3z}}{2z(z-1)}dz\right)=\int_{4/3}^w\textrm{Re}\left(\frac{\sqrt{4-3z}}{2z(z-1)}\right)dx- \textrm{Im}\left(\frac{\sqrt{4-3z}}{2z(z-1)}\right)dy=0$$ which we can do by following the vector field $$\left(\textrm{Im}\left(\frac{\sqrt{4-3z}}{2z(z-1)}\right),\textrm{Re}\left(\frac{\sqrt{4-3z}}{2z(z-1)}\right)\right)$$ from the point $z=4/3$. One part will be the ray $[4/3,\infty)$. The other part will be the interface $\gamma^1$ which will circle around $1$ and will have the shape described in Figure \ref{fig:4}. The interface $\gamma^0$ will come from infinity following the ray $[4/3,\infty)$, then circle clockwise around $D^1$, and then go back towards infinity along the ray $[4/3,\infty)$. 

\subsection{Example with one droplet in $\mathbb{C}\setminus\{1\}$} 

Letting $a_0=1, a_1=0$ and $a_{\infty}=0$ in (\ref{eq:fi}) we get $$\vp(z)=\frac{1}{4z^2(1-z)}dz^2.$$ This quadratic differential corresponds to a stationary solution $(D^0,\gamma^0)$ of the competitive Hele-Shaw problem in $\CC\setminus\{1\}$ with the source divisor $d^0=1\cdot 0$, and it is not hard to see that $D^0=\CC\setminus[1,\infty)$. The corresponding half-translation surface is easily seen to be given by $S^0:=\{ \Re(z)\leq 0, 0\leq \Im(z)\leq \pi\}$ with the top and bottom of each half strip being identified by the translation $z\mapsto z+\pi i$, and the vertical line segment $[0,\pi i/2]$ being identified with $[\pi i/2,\pi i]$ by the minus-translation $z\mapsto -z+\pi i$. Here the point $\pi/2\in S^0$ corresponds to the point $1\in \CC$, while $0,\pi i\in S^0$ both correspond to the point at infinity.

\subsection{Examples with four droplets in $\hat{\mathbb{C}}$}\label{sect:fourdroplets}

Some very special and symmetric examples of stationary solutions with four droplets on $\hat{\mathbb{C}}$ can be found in the following way.

Let $x\in (0,1)$ and denote $$R_x(z):=\log \left|\frac{1-z x}{z-x}\right|.$$ On $\mathbb{D}$ we have that $R_x(z)=G_{\mathbb{D},x}(z)$, and from the observation that $R_x(1/z)=-R_x(z)$ we get that on $\hat{\mathbb{C}}\setminus \overline{\mathbb{D}}$, $R_x(z)=-G_{\hat{\mathbb{C}}\setminus \overline{\mathbb{D}},1/x}(z)$. Let us also note that $R_x(\bar{z})=R_x(z)$.

Now take $x_1,x_2\in (-1,1)$, $a,b\in\mathbb{R}_+$ and consider the function $$R(z):=aR_{x_1}(z)-bR_{x_2}(z).$$ $R$ is then harmonic on $\hat{\mathbb{C}}$ except having $+\infty$ poles at $x_1$ and $1/x_2$ of order $a$ and $b$ respectively, and having $-\infty$ poles at $x_2$ and $1/x_1$ of order $b$ and $a$ respectively. It is also clear that $R=0$ on the unit circle. 

Let us define $D^1:=\{|z|<1, R(z)>0\}, D^2:=\{|z|<1, R(z)<0\}, D^3:=\{|z|>1, R(z)<0\}$ and $D^4:=\{|z|>1, R(z)>0\}$. Since $R(1/z)=-R(z)$ we have that $z\in D^1$ iff $1/z\in D^3$ and similarly $z\in D^2$ iff $1/z\in D^4$. 

If all $D^i$ happen to be simply connected, then by the same argument as in the proof of Proposition \ref{prop:level-line} $(D^i,\partial D^i)$ will be a stationary solution to the competitive Hele-Shaw problem with source divisors $d^1=a\cdot x_1, d^2=b\cdot x_2, d^3=a\cdot 1/x_1$ and $d^4=b\cdot 1/x_2$.   

Note that $R$ is harmonic on the simply connected set $E:=\mathbb{D}\setminus ((-1,x_1]\cup [x_2,1)])$, and by the maximum principle it then follows that both $D^1\cap E$ and $D^2\cap E$ are simply connected. It is however not necessarily true that both $D^1$ and $D^2$ are simply connected. If e.g. $b$
is much smaller than $a$ it can happen that the boundary of $D^2$ does not intersect the boundary of $\partial \mathbb{D}$. But it is easy to see that for $a,b$ fixed, then there exists an $\epsilon>0$ such that if $|x_1+1|<\epsilon$ and $|x_2-1|<\epsilon$ then $R>0$ on $(-1,x_1)$ while $R<0$ on $(x_2,1)$. From this it then follows that both $D^1$ and $D^2$ are simply connected. 

Let us also note that if $D^1$ and $D^2$ are simply connected, then so are $D^3$ and $D^4$.

So at least as long as $|x_1+1|<\epsilon$ and $|x_2-1|<\epsilon$ (where $\epsilon$ depends on the choice of $a,b$) we get a stationary solution to the competitive Hele-Shaw problem with four droplets which is symmetric under inversion $z\mapsto 1/z$ as well as under  reflection in the real axis $z\mapsto \bar{z}$. 

\begin{figure}[h]

\centering

\begin{subfigure}{.8\textwidth}
\centering
    \includegraphics[width=.8\linewidth]{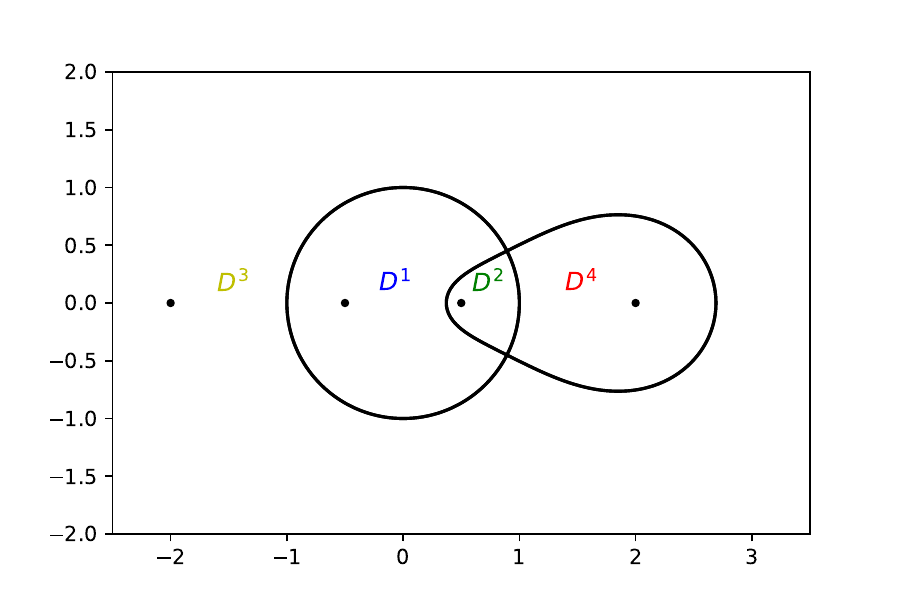}
    \caption{Source points with weight $6$ at $-2$ and $-1/2$ and with weight $1$ at $1/2$ and $2$.}
    \label{fig:135}
\end{subfigure}%

\begin{subfigure}{.8\textwidth}
\centering
    \includegraphics[width=.8\linewidth]{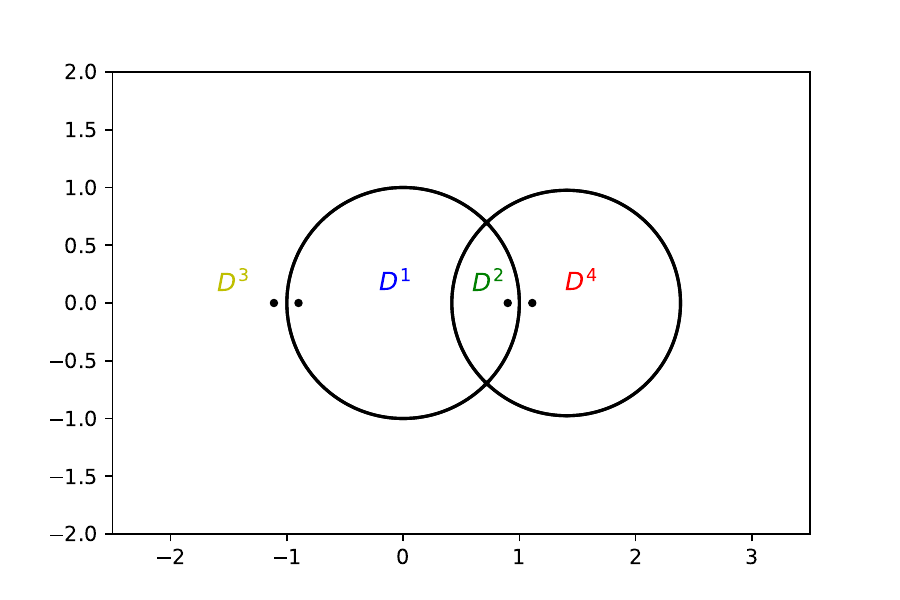}
    \caption{Source points with weight $6$ at $-10/9$ and $-9/10$ with weight $1$ at $9/10$ and $10/9$. The interface separating $D^1$ from $D^2$ is very close to a hyperbolic geodesic in $\mathbb{D}$.}
    \label{fig:246}
\end{subfigure}

\caption{Two examples with four droplets in $\hat{\mathbb{C}}$.}

    \label{fig:113}

\end{figure}

Figure \ref{fig:113}(a) shows the contours of the four droplets when $x_1=-0.5$, $x_2=0.5$, $a=6$ and $b=1$, while in Figure \ref{fig:113}(b) we still have $a=6$ and $b=1$ but now $x_1=-0.9$ and $x_2=0.9$,

If one restricts to the unit disc we only see the droplets $D^1$ and $D^2$, and this corresponds to the situation studied in \cite{GP}.

\section{Lattice models and simulations}\label{sect:discretemodels}
\paragraph{Internal DLA.} Let $K_0 = \{0\}$. The internal DLA (Diffusion Limited Aggregation) cluster $(K_n)_{n \ge 0}$ is a subset of the vertices in $\mathbb{L}$ and grows from $K_n$ to $K_{n+1}$ by adding the first vertex visited by a random walk on $\mathbb{L}$ started from $0$ and stopped when exiting the cluster $K_n$. It can be viewed as a discrete version of classical Hele-Shaw flow with injection from $0$. 

\paragraph{Competitive erosion.}
Competitive erosion \cite{GLPP} is a variation of internal DLA with several different competing clusters. In each starting cluster $C^i_0$ choose a source point $z^i$, and at time one we start a random walk at $z^1$ and let $p$ be the first point visited not contained in $C^1_0$. Then $(C^1_1)':=C^1_0 \cup \{p\}$ while for $i\neq 1$ $(C^i_0)':=C^i_0 \smallsetminus\{p\}$. Next start a random walk at $z^2$ and let $p$ be the first visited point not contained in $(C^2_0)'$. We then let $(C^2_1)':=(C^2_0)'\cup\{p\}$ and $(C^i_0)':=C^i_0 \smallsetminus\{p\}$ for $i\neq 2$. Iterate for all $i$ and in the end let $C^i_1:=(C^i_1)'$. This process is then repeated to get $C^i_k$ for each $k$.

One may have several source points in each cluster, and one can also assign a weight $a\in \NN$ to a given source point $x$, e.g., using Poisson clocks.

\paragraph{Interface erosion on square tiled surfaces.}
Suppose $\Sigma$ is a square tiled surface, i.e., a Riemann surface which can be represented as a translation surface where each polygonal piece is a unit square. This naturally gives rise to a square lattice graph $\mathbb{L}$ on $\Sigma$. Given this the model can be defined analogously to the case of $\Sigma=\mathbb{C}$ (see Definition \ref{df:IE}.) Note that if $\mathbb{L}$ is a square lattice graph on $\Sigma$, then for any $N\in \NN$ we get a refined square lattice graph $\mathbb{L}_N$ on $\Sigma$ simply by dividing each unit square of $\mathbb{L}$ into $N^2$ subsquares.

\paragraph{Long-term behavior.}
For Internal DLA, it is well-known that, with probability one, the growing cluster will contain any bounded subset, given enough time. (However, it is interesting to study the shape of the rescaled cluster, see, e.g., \cite{LBG}.) For competitive erosion the long-term behavior of the clusters is more complicated and the question was considered in \cite{GLPP} and \cite{GP}. 

In \cite{GP} Ganguly-Peres consider competitive erosion on a $\frac{1}{N}\mathbb{Z}^n$ lattice approximation of a smooth simply connected domain with two competing clusters with source points close to two given boundary points $z$ and $w$. It is shown that in the small-mesh limit the interface between the two clusters approaches a level set of the Green's function $G$ with Neumann boundary conditions on $D$ with source and sink at $z$ and $w$ respectively. The limiting interface is a hyperbolic geodesic separating $z$ and $w$. While this setting is somewhat different from ours (for instance, the surface $\mathbb{D}$ has a boundary), the limiting configuration can be obtained from a quadratic differential after reflection, see Section~\ref{sect:fourdroplets}. 

We have the following conjecture for interface erosion.
\begin{conj} \label{conj:comper}
    Let $\Sigma$ be a square tiled surface with a square lattice graph $\mathbb{L}$ and let $g$ be the induced metric on $\Sigma$. Assume that for all $N\in \NN$, $(D_t^{i,N},\gamma_t^{i,N})$ is the result of running interface erosion on the square lattice graph $\mathbb{L}_N$ and source divisors $d^i$ starting from $(D_0^{i,N},\gamma_0^{i,N})$. Furthermore assume that for each $i$ the droplet interfaces $\gamma_0^{i,N}$ converge to a given droplet interface $\gamma^i_0$ as $N\to \infty$. Then as $N\to \infty$, $(D_{N^2t}^{i,N},\gamma_{N^2t}^{i,N})$ converge in probability towards a solution $(D^i_t,\gamma^i_t)$ to the competitive Hele-Shaw problem on $(\Sigma,g)$ with source divisors $d^i$.
\end{conj}

In particular, if the lattice is sufficiently fine we expect that for large $N$ \emph{and} $t$ the discrete droplet interfaces will be close to a stationary solution of the associated competitive Hele-Shaw problem.

Our conjecture fits well with the above mentioned result of Ganguly-Peres, as is explained in Section~\ref{sect:fourdroplets}.

\end{document}